# TWISTED-PRODUCT CATEGORICAL BUNDLES

SAIKAT CHATTERJEE, AMITABHA LAHIRI, AND AMBAR N. SENGUPTA

Abstract. Categorical bundles provide a natural framework for gauge theories involving multiple gauge groups. Unlike the case of traditional bundles there are distinct notions of triviality, and hence also of local triviality, for categorical bundles. We study categorical principal bundles that are product bundles in the categorical sense, developing the relationship between functorial sections of such bundles and trivializations. We construct functorial cocycles with values in categorical groups using a suitable family of locally defined functors on the object space of the base category. Categorical product bundles being too rigid to give a widely applicable model for local triviality, we introduce the notion of a twisted-product categorical bundle. We relate such bundles to decorated categorical bundles that contain more information, specifically parallel transport data.

## 1. Introduction

The purpose of this paper is to explore the nature of product categorical bundles and propose the notion of *twisted-product* categorical bundles, a richer and more widely applicable notion that could be used to formulate local triviality for categorical principal bundles. The main motivation for this exploration arises from there being more than one useful notion of triviality, and hence also of local triviality, for categorical bundles. For this reason the definition of categorical principal bundle we work with does not include a local triviality requirement. In these introductory paragraphs we present a discussion of our investigations phrased in broadly understandable terms, leaving exact definitions to later sections.

The language of category theory has proved very useful in the formulation of gauge theories that involve multiple gauge groups, each such gauge group describing the interaction of gauge fields with particles or higher-dimensional entities. Thus such "higher gauge theories" inspire a categorical formulation of the notions of principal bundles and parallel transport. In the simplest case we would view the points of a manifold as objects of a category and paths would be the morphisms (this is a naive first view). At the next higher level, paths would be the objects and morphisms would be paths of paths. Thus a traditional bundle has as counterpart a projection functor

$$\pi : \mathbf{P} \longrightarrow \mathbf{B},$$

where $\mathbf{P}$ and $\mathbf{B}$ are categories that are counterparts of the bundle space and the base space, respectively. The geometry of a categorical principal bundle also includes a *categorical group* in place of the traditional structure group; in a categorical group there is both a group comprised of objects and a group comprised of morphisms, with these structures intertwined in a consistent way. A principal categorical bundle is then given by a surjective functor $\mathbf{P} \longrightarrow \mathbf{B}$ along with a right action (in a suitable sense) of the categorical group $\mathbf{G}$ on $\mathbf{P}$, satisfying conditions analogous to those satisfied by traditional principal bundles. In section 2 we state the basic definitions for categorical principal bundles, along with a review of other related notions.

As noted above, thinking of the morphisms of $\mathbf{P}$ simply as arbitrary paths on $P$ is a naive first view. There is far too much freedom and to little structure in such a large space for our objectives. In the examples of most interest to us a morphism of $\mathbf{P}$ corresponds to a path $\overline{\gamma}$ on $P$ that is *horizontal* with respect to some connection $\overline{A}$ on $P$ (thus describing parallel-transport along some path on $B$), decorated by an element $h$ drawn from the second gauge group $H$. Figure 4 illustrates such a morphism $(\overline{\gamma}, h)$. The geometric information, coming from parallel transport, is encoded in the categorical structure in the target of this morphism and therefore in the way composition is defined.

The simplest example of a principal categorical comes from the categorical product $\mathbf{U} \times \mathbf{G}$; it is given by the projection onto the base category $\mathbf{U}$:

$$\mathbf{U} \times \mathbf{G} \longrightarrow \mathbf{U},$$

along with the obvious action of $\mathbf{G}$ on $\mathbf{U}\times\mathbf{G}$. In classical bundle theory, the product bundle, despite its bare and simple nature, is important in that it serves as a local model for more general *locally trivial* bundles. For categorical principal bundles the notion of local triviality is less obvious. A major motivation for our investigations in this paper lies in providing such a local model for categorical principal bundles.

In sections 3 and 4 we explore categorical product bundles and prove results that are the categorical counterparts of results on traditional product topological bundles. As we have noted above, the motivation is not so much the focus on product bundles as the understanding of what happens locally in more general bundles. For example, traditionally, sections of a product bundle $U \times G \longrightarrow U$, where $G$ is a group, correspond to mappings $U \longrightarrow G$. We show in Proposition 3.4 that in the categorical setting the collection of sections forms a categorical group whose objects are functors $\mathbf{U} \longrightarrow \mathbf{G}$ and whose morphisms are natural transformations between such functors. The proof involves all aspects of the structure categorical group $\mathbf{G}$.

A traditional principal bundle is given by a smooth projection map

$$\pi : P \longrightarrow B,$$

along with a free right action of a Lie group $G$ on $P$, acting transitively on each fiber $\pi^{-1}(b)$. Moreover, such a bundle is locally trivial in the sense that there is an open covering $\{U_\alpha\}_{\alpha \in I}$ of $B$ such that the bundle viewed over any $U_\alpha$ is trivial; this means that there is a diffeomorphism

$$\phi_\alpha : U_\alpha \times G \longrightarrow \pi^{-1}(U_\alpha) \tag{1.1}$$

that respects the action of $G$ and maps the fiber $\{x\} \times G$ onto the fiber $\pi^{-1}(x)$ for each point $x \in U_\alpha$. Switching from $U_\alpha$ to $U_\beta$, which intersects $U_\alpha$, is described then by a mapping

$$g_{\alpha\beta} : U_\alpha \cap U_\beta \longrightarrow G,$$

specified by

$$\phi_\alpha^{-1}\Big(\phi_\beta(x, e)\Big) = (x, e)g_{\alpha\beta}(x), \tag{1.2}$$

for all $x \in U_\alpha \cap U_\beta$. The system of functions $\{g_{\alpha\beta}\}$ is a *cocycle*, satisfying the condition

$$g_{\alpha\beta}(x)g_{\beta\gamma}(x) = g_{\alpha\gamma}(x) \tag{1.3}$$

for all $x \in U_\alpha \cap U_\beta \cap U_\gamma$. Section 5 is devoted to studying a functorial counterpart of this notion of cocycle.

We turn next, in section 6, to the study of a type of categorical principal bundle that is richer in structure than product categorical bundles. For a *twisted-product categorical bundle*

$$\mathbf{U} \times_\eta \mathbf{G} \longrightarrow \mathbf{U},$$

the object and morphism sets of $\mathbf{U} \times_\eta \mathbf{G}$ are the same as for the product category $\mathbf{U} \times \mathbf{G}$; however, the target maps and composition are defined differently. There is more structure; the idea is that the 'twist' $\eta$ encodes a connection form on a traditional product bundle $U \times G \longrightarrow U$.

Let us note that *we do not use the term 'twist' to mean topological nontriviality* of a bundle; it is, instead, a more subtle notion that is expressed, for instance, in the way compositions of morphisms of $\mathbf{P}$ is defined. Our objective here is to create a concept that is a useful categorical counterpart of product bundles to be used for local triviality of categorical principal bundles.

There is too extensive a literature in category theoretic geometry to offer a comprehensive review here. Broadly related to our work are those of Abbaspour and Wagemann [1], Attal [3,4], Baez et al. [5,6], Barrett [7], Bartels [8], Breen and Messing [9], Gawędzki and Reis [16], Picken et al. [10,19,20], Parzygnat [22], Schreiber and Waldorf [23,24], Soncini and Zucchini [25], and Viennot [26].

## 2. Categorical principal bundles

In this section we review the fundamentals of the theory of categorical principal bundles and categorical connections in the framework developed in our work [12]. Let us note once and for all that for our categories, objects and morphisms form sets; thus, associated to a category $\mathbf{C}$ there is the object set $\mathrm{Obj}(\mathbf{C})$ and the morphism set $\mathrm{Mor}(\mathbf{C})$.

### 2.1. Categorical groups.

We begin with a summary of the essentials of categorical groups. The concepts and results we discuss here are standard. References that may be consulted for more information include [15,17,18].

By a *categorical group* $\mathbf{G}$ we mean a category $\mathbf{G}$ along with a functor

$$\mathbf{G} \times \mathbf{G} \longrightarrow \mathbf{G} \tag{2.1}$$

that makes both $\mathrm{Obj}(\mathbf{G})$ and $\mathrm{Mor}(\mathbf{G})$ groups. Useful consequences include the following: (i) the identity-assigning map

$$\mathrm{Obj}(\mathbf{G}) \longrightarrow \mathrm{Mor}(\mathbf{G}) : x \mapsto 1_x$$

is a homomorphism, (ii) the source and target maps

$$s, t : \mathrm{Mor}(\mathbf{G}) \longrightarrow \mathrm{Obj}(\mathbf{G}) \tag{2.2}$$

are both homomorphisms; (iii) the following exchange law, expressing functoriality of (2.1), holds

$$(\phi_2\psi_2) \circ (\phi_1\psi_1) = (\phi_2 \circ \phi_1)(\psi_2 \circ \psi_1) \tag{2.3}$$

for all $\phi_1, \phi_2, \psi_1, \psi_2 \in \mathrm{Mor}(\mathbf{G})$.

A *categorical Lie group* is a categorical group $\mathbf{G}$ for which both $\mathrm{Obj}(\mathbf{G})$ and $\mathrm{Mor}(\mathbf{G})$ are Lie groups and the maps $s$, $t$ and $x \mapsto 1_x$ are smooth.

A *crossed module* is a pair of groups $G$ and $H$, along with maps

$$\alpha : G \times H \longrightarrow H : (g, h) \mapsto \alpha_g(h) \quad \text{and} \quad \tau : H \longrightarrow G,$$

where $\tau$ is a homomorphism, $\alpha_g$ is an automorphism of $H$ for each $g \in G$, and the map $g \mapsto \alpha_g \in \mathrm{Aut}(H)$ is a homomorphism. The target map $\tau$ and the map $\alpha$ interact through the relation

$$\alpha_{\tau(h)}(h') = hh'h^{-1} \qquad \text{for all } h, h' \in H. \tag{2.4}$$

We can also write this as

$$\tau(h)h'\tau(h)^{-1} = hh'h^{-1} \qquad \text{for all } h, h' \in H. \tag{2.5}$$

When $G$ and $H$ are Lie groups, and $\alpha$ and $\tau$ are smooth, $(G, H, \alpha, \tau)$ is called a *Lie crossed module.*

There is a bijective correspondence between categorical (Lie) groups and (Lie) crossed modules. Given a categorical group $\mathbf{G}$, we take $G = \mathrm{Obj}(\mathbf{G})$, $H = \ker s$, $\tau = t|H$, and

$$\alpha_g(h) = 1_g h 1_g^{-1}$$

for all $g \in G$ and $h \in H$. Then $\mathrm{Mor}(\mathbf{G})$ is isomorphic to the semidirect product $H \rtimes_\alpha G$ via the map

$$\mathrm{Mor}(\mathbf{G}) \longrightarrow H \rtimes_\alpha G : \phi \mapsto \big(\phi 1_{s(\phi)^{-1}}, s(\phi)\big). \tag{2.6}$$

The target map $t$, viewed as a mapping $H \rtimes_\alpha G \longrightarrow G$, is given by

$$t(h, g) = \tau(h)g \qquad \text{for all } (h, g) \in H \rtimes_\alpha G. \tag{2.7}$$

Let us also note that the binary operation in the semidirect product $H \rtimes_\alpha G$ is given by

$$(h_2, g_2)(h_1, g_1) = \big(h_2\alpha_{g_2}(h_1), g_2g_1\big) \tag{2.8}$$

for all $(h_2, g_2), (h_1, g_1) \in H \rtimes_\alpha G$.

It will be convenient to identify $h \in H$ with $(h, e) \in H \rtimes_\alpha G$ and $g \in G$ with $(e, g) \in H \rtimes_\alpha G$; then

$$hg = (h, g) \tag{2.9}$$

and the product (2.8) becomes

$$h_2g_2h_1g_1 = h_2\alpha_{g_2}(h_1)g_2g_1. \tag{2.10}$$

The automorphism $\alpha_g : H \longrightarrow H$ is then just conjugation:

$$\alpha_g(h) = ghg^{-1} \qquad \text{for all } h \in H \text{ and } g \in G. \tag{2.11}$$

From the product formula in (2.10) we see then that

$$\text{the } H\text{-component of } (h_2, g_2)(h_1, g_1) \text{ is } h_2g_2h_1g_2^{-1}. \tag{2.12}$$

The composition of morphisms in $\mathrm{Mor}(\mathbf{G}) \simeq H \rtimes_\alpha G$ is given by

$$(h_2, g_2) \circ (h_1, g_1) = (h_2h_1, g_1); \tag{2.13}$$

this composition is defined only for those $(h_2, g_2), (h_1, g_2) \in H \rtimes_\alpha G$ for which

$$\tau(h_1)g_1 = t(h_1, g_1) = s(h_2, g_2) = g_2.$$

2.2. **Categorical principal bundles.** A *right action* of a categorical group $\mathbf{G}$ on a category $\mathbf{P}$ is a functor

$$\mathbf{P} \times \mathbf{G} \longrightarrow \mathbf{P}$$

which at the object level gives a right action of $\mathrm{Obj}(\mathbf{G})$ on $\mathrm{Obj}(\mathbf{P})$ and on the morphism level gives a right action of $\mathrm{Mor}(\mathbf{G})$ on $\mathrm{Mor}(\mathbf{P})$. We say that this action is *free* if it is free both at the object level and at the morphism level.

A *categorical principal bundle* comprises categories $\mathbf{P}$ and $\mathbf{B}$ along with a right action of a categorical group $\mathbf{G}$ on $\mathbf{P}$ satisfying the following conditions:

(b1) $\pi$ is surjective both at the level of objects and at the level of morphisms;
(b2) the action of $\mathbf{G}$ on $\mathbf{P}$ is free on both objects and morphisms;
(b3) the action of $\mathrm{Obj}(\mathbf{G})$ preserves the fiber $\pi^{-1}(b)$ and its action on this fiber is transitive, for each object $b \in \mathrm{Obj}(\mathbf{B})$; the action of $\mathrm{Mor}(\mathbf{G})$ preserves the fiber $\pi^{-1}(\phi)$ and its action on this fiber is transitive, for each morphism $\phi \in \mathrm{Mor}(\mathbf{B})$.

(This definition is from [12].)

In the central examples of interest both $\mathrm{Obj}(\mathbf{P})$ and $\mathrm{Obj}(\mathbf{B})$ are manifolds, while the morphisms arise from certain types of paths on these manifolds. (The manifold itself might be a space of paths). In this context, of course, the functors all arise from continuous maps. However, we do not explicitly use any topological results and therefore have not spelled out a topology in this work.

A categorical principal bundle contains much more than just the projection functor: the action of the categorical group $\mathbf{G}$ on $\mathbf{P}$ includes a considerable amount of structure because it operates on both morphisms and objects in a way that is consistent with the source and target maps. In examples of interest the source and target maps themselves contain nontrivial data such as information on parallel transport arising from a connection. This is discussed in more detail in section 6.3. Unlike the case of classical bundles there are *multiple distinct notions of triviality*; for this reason we do not build a requirement of local triviality into the definition of a categorical principal bundle.

## 3. The product categorical bundle

Our definition of categorical principal bundle does not require local triviality. In order to formulate a notion of local triviality we need to study a categorical counterpart of a trivial bundle. Traditionally a principal bundle is trivial if it is isomorphic, as a principal bundle, to a product bundle. Moreover, in the traditional setting, a principal bundle is trivial if and only if it has a global section (of the relevant degree of smoothness). In this section we will study the categorical analog of a product bundle and related ideas and results.

A *product categorical principal bundle* is given by the projection on the first factor

$$\mathbf{U} \times \mathbf{G} \longrightarrow \mathbf{U},$$

where $\mathbf{U}$ is any category and $\mathbf{G}$ is a categorical group acting on the right on $\mathbf{U} \times \mathbf{G}$ in the obvious way.

An object of $\mathbf{U} \times \mathbf{G}$ is of the form

$$(a, g) \in U \times G = \mathrm{Obj}(\mathbf{U}) \times \mathrm{Obj}(\mathbf{G}).$$

A morphism of $\mathbf{U} \times \mathbf{G}$ is of the form

$$(\gamma, h, g) \in \mathrm{Mor}(\mathbf{U}) \times (H \rtimes_\alpha G).$$

Its source and target are given by

$$s(\gamma, h, g) = \big(s(\gamma), g\big) \qquad \text{and} \qquad t(\gamma, h, g) = \big(t(\gamma), \tau(h)g\big). \tag{3.1}$$

The right action at the level of objects is given by

$$(a, g)g_1 = (a, gg_1)$$

and, recalling (2.8), at the level of morphisms by

$$(\gamma, h, g)(h_1, g_1) = \big(\gamma, h\alpha_g(h_1), gg_1\big), \tag{3.2}$$

which we can also display as

$$(\gamma, hg)h_1g_1 = (\gamma, hgh_1g_1). \tag{3.3}$$

Composition of morphisms is given by composition in each component:

$$(\gamma_2, h_2, g_2) \circ (\gamma_1, h_1, g_1) = (\gamma_2 \circ \gamma_1, h_2h_1, g_1), \tag{3.4}$$

where we have used the composition law in $\mathbf{G}$, given by (2.13), in the second component, with $g_2 = \tau(h_1)g_1$.

A *section* of a traditional product bundle

$$\pi : U \times G \longrightarrow U : (a, g) \mapsto a$$

is a mapping

$$s : U \longrightarrow U \times G$$

for which $\pi \circ s(a) = a$ for all $a \in U$. Thus a section is of the form

$$U \longrightarrow U \times G : a \mapsto \big(a, \psi(a)\big),$$

where $\psi : U \longrightarrow G$ is a mapping. The set of all mappings $U \longrightarrow G$ is a group under pointwise multiplication.

Analogously, for the categorical product bundle

$$\mathbf{U} \times \mathbf{G} \longrightarrow \mathbf{U}$$

we can view sections as corresponding to functors

$$\mathbf{U} \longrightarrow \mathbf{G}.$$

We will study the nature of such functors and show that *these functors along with natural transformations between them form a categorical group.*

We discuss our results first, collecting the proofs all together after this discussion.

Our first step is to specify a concrete way of encoding functors $\mathbf{U} \longrightarrow \mathbf{G}$. As we see in the following proposition, such a functor $\mathbf{\Psi}$ can be understood completely in terms of a function $g : \mathrm{Obj}(\mathbf{U}) \longrightarrow G$ and a function $h : \mathrm{Mor}(\mathbf{U}) \longrightarrow H$.

**Proposition 3.1.** *Let* $\mathbf{U}$ *be any category and* $\mathbf{G}$ *a categorical group. Let* $U = \mathrm{Obj}(\mathbf{U})$, $G = \mathrm{Obj}(\mathbf{G})$ *and* $H = \ker s : \mathrm{Mor}(\mathbf{G}) \longrightarrow G$, *as usual, so that* $\mathrm{Mor}(\mathbf{G}) = H \rtimes_\alpha G$, *for a homomorphism* $\alpha : G \longrightarrow \mathrm{Aut}(H)$. *For notational convenience we write the target and source of a morphism as*

$$\gamma_1 = t(\gamma) \text{ and } \gamma_0 = s(\gamma).$$

*Suppose*

$$\mathbf{\Psi} : \mathbf{U} \longrightarrow \mathbf{G}$$

*is a functor. Then there are functions*

$$g : U \longrightarrow G \qquad \text{and} \qquad h_{\mathbf{\Psi}} : \mathrm{Mor}(\mathbf{U}) \longrightarrow H \tag{3.5}$$

*such that*

$$\begin{aligned} \mathbf{\Psi}(a) &= g(a) \qquad \text{for all } a \in \mathrm{Obj}(\mathbf{U}) \\ \mathbf{\Psi}(\gamma) &= \big(h_{\mathbf{\Psi}}(\gamma), g(\gamma_0)\big) \qquad \text{for all } \gamma \in \mathrm{Mor}(\mathbf{U}), \end{aligned} \tag{3.6}$$

*wherein* $\gamma_0 = s(\gamma)$; *moreover, the function* $h_{\mathbf{\Psi}}$ *has the following properties:*

$$h_{\mathbf{\Psi}}(\gamma' \circ \gamma) = h_{\mathbf{\Psi}}(\gamma') h_{\mathbf{\Psi}}(\gamma) \tag{3.7}$$

*for all* $\gamma, \gamma' \in \mathrm{Mor}(\mathbf{U})$ *for which the composite* $\gamma' \circ \gamma$ *is defined, and*

$$\tau\Big(h_{\mathbf{\Psi}}(\gamma)\Big) = g\big(\gamma_1\big) g\big(\gamma_0\big)^{-1}. \tag{3.8}$$

*In the converse direction, suppose* $h : U \longrightarrow H$ *is a function and* $g : U \longrightarrow G$ *the function given by*

$$g = \tau \circ h : U \longrightarrow G.$$

*Then the assigments:*

$$\begin{aligned} &\mathbf{\Psi}(a) = g(a) \qquad \textit{for all } a \in U; \\ &\mathbf{\Psi}(\gamma) = \big(h(\gamma_1)h(\gamma_0)^{-1}, g(\gamma_0)\big) \qquad \textit{for all } \gamma \in \mathrm{Mor}(\mathbf{U}) \end{aligned} \tag{3.9}$$

*specify a functor* $\mathbf{\Psi} : \mathbf{U} \longrightarrow \mathbf{G}$*.*

Our next observation is that the pointwise-defined product of two functors $\mathbf{U} \longrightarrow \mathbf{G}$ is another such functor, and indeed the such functors form a group under this pointwise multiplication.

**Proposition 3.2.** *Suppose*

$$\mathbf{\Psi}_1, \mathbf{\Psi}_2 : \mathbf{U} \longrightarrow \mathbf{G}$$

*are functors, where* $\mathbf{U}$ *is any category and* $\mathbf{G}$ *is a categorical group. Then the pointwise product*

$$\mathbf{\Psi}_2\mathbf{\Psi}_1 : \mathbf{U} \longrightarrow \mathbf{G}$$

*is also a functor. The functors* $\mathbf{U} \longrightarrow \mathbf{G}$ *form a group under pointwise multiplication.*

A *natural transformation* $\mathbb{T}$ from a functor $\mathbf{\Psi}_1 : \mathbf{U} \Rightarrow \mathbf{G}$ to a functor $\mathbf{\Psi}_2 : \mathbf{U} \longrightarrow \mathbf{G}$ is, by definition, an assignment to each object $a \in \mathrm{Obj}(\mathbf{U})$ a morphism

$$\mathbb{T}(a) : \mathbf{\Psi}_1(a) \longrightarrow \mathbf{\Psi}_2(a)$$

such that the diagram

$$\begin{array}{ccc} \mathbf{\Psi}_1(a) & \xrightarrow{\mathbf{\Psi}_1(f)} & \mathbf{\Psi}_1(b) \\ \Big\downarrow{\scriptstyle \mathbb{T}(a)} & & {\scriptstyle \mathbb{T}(b)}\Big\downarrow \\ \mathbf{\Psi}_2(a) & \xrightarrow[\mathbf{\Psi}_2(f)]{} & \mathbf{\Psi}_2(b) \end{array} \tag{3.10}$$

commutes for all morphisms $f : a \longrightarrow b$ in $\mathrm{Mor}(\mathbf{U})$.

The following result gives an explicit description of such transformations; every such transformation arises from a function $\mathrm{Obj}(\mathbf{U}) \longrightarrow H$.

**Proposition 3.3.** *Suppose*

$$\mathbf{\Psi}_1, \mathbf{\Psi}_2 : \mathbf{U} \longrightarrow \mathbf{G}$$

*are functors, where* $\mathbf{U}$ *is a category and* $\mathbf{G}$ *is a categorical group, and suppose*

$$\mathbb{T} : \mathbf{\Psi}_1 \Rightarrow \mathbf{\Psi}_2$$

*is a natural transformation. Let $(G, H, \alpha, \tau)$ be the crossed module associated to $\mathbf{G}$, so that $\mathrm{Mor}(\mathbf{G}) = H \rtimes_\alpha G$ and $\mathrm{Obj}(\mathbf{G}) = G$. Then there is a function*

$$h_{\mathbb{T}} : \mathrm{Obj}(\mathbf{U}) \longrightarrow H$$

*such that*

$$\mathbf{\Psi}_2(a) = \tau\big(h_{\mathbb{T}}(a)\big)\mathbf{\Psi}_1(a) \tag{3.11}$$

*and*

$$h_{\mathbf{\Psi}_2}(f) = h_{\mathbb{T}}(b) h_{\mathbf{\Psi}_1}(f) h_{\mathbb{T}}(a)^{-1} \tag{3.12}$$

*for all morphisms $f : a \longrightarrow b$ in $\mathrm{Mor}(\mathbf{U})$, where the notation $h_{\mathbf{\Psi}_i}(f)$ is specified by denoting the morphism $\mathbf{\Psi}_i(f) \in \mathrm{Mor}(\mathbf{G}) = H \rtimes_\alpha G$ as $\big(h_{\mathbf{\Psi}_i}(f), g_{\mathbf{\Psi}_i}(f)\big)$:*

$$\mathbf{\Psi}_i(f) = \big(h_{\mathbf{\Psi}_i}(f), g_{\mathbf{\Psi}_i}(f)\big).$$

Now let us consider natural transformations $\mathbb{T} : \mathbf{\Psi}_1 \Rightarrow \mathbf{\Psi}_2$ and $\mathbb{T}' : \mathbf{\Psi}'_1 \Rightarrow \mathbf{\Psi}'_2$ between functors $\mathbf{U} \longrightarrow \mathbf{G}$. Then we can define a pointwise product $\mathbb{T}'\mathbb{T}$ that associates to each object $a \in \mathrm{Obj}(\mathbf{U})$ the morphism

$$(\mathbb{T}'\mathbb{T})(a) \stackrel{\text{def}}{=} \mathbb{T}'(a)\mathbb{T}(a) : \mathbf{\Psi}'_1(a)\mathbf{\Psi}_1(a) \longrightarrow \mathbf{\Psi}'_2(a)\mathbf{\Psi}_2(a).$$

From the two commuting diagrams in $\mathrm{Mor}(\mathbf{G})$

$$\begin{array}{ccc} \mathbf{\Psi}_1(a) & \xrightarrow{\mathbf{\Psi}_1(f)} & \mathbf{\Psi}_1(b) \\ \downarrow{\scriptstyle \mathbb{T}(a)} & & \downarrow{\scriptstyle \mathbb{T}(b)} \\ \mathbf{\Psi}_2(a) & \xrightarrow[\mathbf{\Psi}_2(f)]{} & \mathbf{\Psi}_2(b) \end{array} \qquad \begin{array}{ccc} \mathbf{\Psi}'_1(a) & \xrightarrow{\mathbf{\Psi}'_1(f)} & \mathbf{\Psi}'_1(b) \\ \downarrow{\scriptstyle \mathbb{T}'(a)} & & \downarrow{\scriptstyle \mathbb{T}'(b)} \\ \mathbf{\Psi}'_2(a) & \xrightarrow[\mathbf{\Psi}'_2(f)]{} & \mathbf{\Psi}'_2(b) \end{array} \tag{3.13}$$

we obtain the 'pointwise product diagram'

$$\begin{array}{ccc} \mathbf{\Psi}_1(a)\mathbf{\Psi}'_1(a) & \xrightarrow{\mathbf{\Psi}_1(f)\mathbf{\Psi}'_1(f)} & \mathbf{\Psi}_1(b)\mathbf{\Psi}'_1(b) \\ \downarrow{\scriptstyle \mathbb{T}'(a)\mathbb{T}(a)} & & \downarrow{\scriptstyle \mathbb{T}'(b)\mathbb{T}(b)} \\ \mathbf{\Psi}_2(a)\mathbf{\Psi}'_2(a) & \xrightarrow[\mathbf{\Psi}_2(f)\mathbf{\Psi}'_2(f)]{} & \mathbf{\Psi}_2(b)\mathbf{\Psi}'_2(b) \end{array}. \tag{3.14}$$

Functoriality of the group operation $\mathbf{G} \times \mathbf{G} \longrightarrow \mathbf{G}$ implies that this diagram commutes.

Thus we have a product operation on the natural transformations between functors $\mathbf{U} \longrightarrow \mathbf{G}$.

If

$$\mathbb{T}_1 : \mathbf{\Psi}_1 \Rightarrow \mathbf{\Psi}_2 \qquad \text{and} \qquad \mathbb{T}_2 : \mathbf{\Psi}_2 \Rightarrow \mathbf{\Psi}_3$$

are natural transformations then the composite

$$\mathbb{T}_2 \circ \mathbb{T}_1 : \mathbf{\Psi}_1 \Rightarrow \mathbf{\Psi}_3$$

is defined by

$$(\mathbb{T}_2 \circ \mathbb{T}_1)(a) = \mathbb{T}_2(a) \circ \mathbb{T}_1(a) : \boldsymbol{\Psi}_1(a) \longrightarrow \boldsymbol{\Psi}_3(a) \qquad \text{for all } a \in \mathrm{Obj}(\mathbf{U}). \tag{3.15}$$

With the composition law (3.15) we obtain a category

$$\mathbf{G}^{\mathbf{U}} \tag{3.16}$$

where objects are functors $\mathbf{U} \longrightarrow \mathbf{G}$ and morphisms are natural transformations between such functors. As we have already seen, both $\mathrm{Mor}(\mathbf{G}^{\mathbf{U}})$ and $\mathrm{Obj}(\mathbf{G}^{\mathbf{U}})$ are groups.

If $\mathbb{T}_1, \mathbb{T}_2$ and $\mathbb{T}'_1, \mathbb{T}'_2$ are natural transformations for which the composites $\mathbb{T}_2 \circ \mathbb{T}_1$ and $\mathbb{T}'_2 \circ \mathbb{T}'_1$ are defined then

$$(\mathbb{T}'_2\mathbb{T}_2) \circ (\mathbb{T}'_1\mathbb{T}_1) = (\mathbb{T}'_2 \circ \mathbb{T}'_1)(\mathbb{T}_2 \circ \mathbb{T}_1). \tag{3.17}$$

We prove this by applying to any object $a \in \mathrm{Obj}(\mathbf{U})$:

$$\begin{aligned}
[(\mathbb{T}'_2\mathbb{T}_2) \circ (\mathbb{T}'_1\mathbb{T}_1)](a) &= (\mathbb{T}'_2\mathbb{T}_2)(a) \circ (\mathbb{T}'_1\mathbb{T}_1)(a) \\
&= \big(\mathbb{T}'_2(a)\mathbb{T}_2(a)\big) \circ \big(\mathbb{T}'_1(a)\mathbb{T}_1(a)\big) \\
&= \big(\mathbb{T}'_2(a) \circ \mathbb{T}'_1(a)\big)\big(\mathbb{T}_2(a) \circ \mathbb{T}_1(a)\big) \\
&\quad \text{(because the group operation on } \mathbf{G} \text{ is functorial)} \\
&= [(\mathbb{T}'_2 \circ \mathbb{T}'_1)(\mathbb{T}_2 \circ \mathbb{T}_1)](a).
\end{aligned} \tag{3.18}$$

It is well known [18, page 42] that there are two compositions, satisfying the exchange law (3.17), for certain categories whose objects are morphisms and whose arrows are natural transformations. However, our situation here is different and the second operation between natural transformations is defined using the group structure on $\mathbf{G}$.

The exchange law (3.17) says that the group operation

$$\mathbf{G}^{\mathbf{U}} \times \mathbf{G}^{\mathbf{U}} \longrightarrow \mathbf{G}^{\mathbf{U}}$$

is functorial.

Thus we have proved:

**Proposition 3.4.** $\mathbf{G}^{\mathbf{U}}$ *is a categorical group.*

We now present the proofs of the results of this section.

*Proof of Proposition 3.1.* We define the function $g$ by the first equation in (3.6). Then, writing $\boldsymbol{\Psi}(\gamma)$ as

$$\boldsymbol{\Psi}(\gamma) = \Big(h_{\boldsymbol{\Psi}}(\gamma), g_{\boldsymbol{\Psi}}(\gamma)\Big), \tag{3.19}$$

we have

$$\begin{aligned} s\big(h_{\mathbf{\Psi}}(\gamma), g_{\mathbf{\Psi}}(\gamma)\big) &= s\mathbf{\Psi}(\gamma) = \mathbf{\Psi}s(\gamma) = g\big(\gamma_0\big) \\ t\big(h_{\mathbf{\Psi}}(\gamma), g_{\mathbf{\Psi}}(\gamma)\big) &= t\mathbf{\Psi}(\gamma) = \mathbf{\Psi}t(\gamma) = g(\gamma_1). \end{aligned} \tag{3.20}$$

From the first equation we have

$$g_{\mathbf{\Psi}}(\gamma) = g(\gamma_0). \tag{3.21}$$

The second equation says

$$\tau\big(h_{\mathbf{\Psi}}(\gamma)\big)g_{\mathbf{\Psi}}(\gamma) = g\big(\gamma_1\big). \tag{3.22}$$

Using the value $g_{\mathbf{\Psi}}(\gamma) = g(\gamma_0)$ we then obtain

$$\tau\big(h_{\mathbf{\Psi}}(\gamma)\big) = g\big(\gamma_1\big)g\big(\gamma_0\big)^{-1}. \tag{3.23}$$

Since $\mathbf{\Psi}$ is a functor we also have

$$\mathbf{\Psi}(\gamma' \circ \gamma) = \mathbf{\Psi}(\gamma') \circ \mathbf{\Psi}(\gamma),$$

whenever the composite $\gamma' \circ \gamma$ is defined in $\mathrm{Mor}(\mathbf{U})$ and so, focusing on what this says about the $H$-component, we see that

$$h_{\mathbf{\Psi}}(\gamma' \circ \gamma) = h_{\mathbf{\Psi}}(\gamma')h_{\mathbf{\Psi}}(\gamma). \tag{3.24}$$

Now we verify the converse implication. Suppose $\mathbf{\Psi} : \mathbf{U} \longrightarrow \mathbf{G}$ is given on objects and morphisms by (3.9). It is clear that $\mathbf{\Psi}$ carries identity morphisms to identity morphisms. Next we verify that $\mathbf{\Psi}$ respects sources and targets:

$$\begin{aligned} s\mathbf{\Psi}(\gamma) &= g(\gamma_0) = \mathbf{\Psi}(s\gamma) \\ t\mathbf{\Psi}(\gamma) &= \tau\Big(h(\gamma_1)h(\gamma_0)^{-1}\Big)g(\gamma_0) \\ &= g(\gamma_1) \\ &= \mathbf{\Psi}(t\gamma). \end{aligned} \tag{3.25}$$

Finally, we verify that $\mathbf{\Psi}$ preserves compositions: if $\gamma, \gamma' \in \mathrm{Mor}(\mathbf{U})$ are such that $t(\gamma) = s(\gamma')$ then

$$\begin{aligned} \mathbf{\Psi}(\gamma') \circ \mathbf{\Psi}(\gamma) &= \Big(h(\gamma_1')h(\gamma_0')^{-1}, g(\gamma_0')\Big) \circ \Big(h(\gamma_1)h(\gamma_0)^{-1}, g(\gamma_0)\Big) \\ &= \Big(h(\gamma_1')h(\gamma_0)^{-1}, g(\gamma_0)\Big) \qquad \text{using } h(\gamma_0') = h(\gamma_1) \\ &= \mathbf{\Psi}(\gamma' \circ \gamma). \end{aligned} \tag{3.26}$$

Thus $\mathbf{\Psi}$ is a functor. $\square$

*Proof of Proposition 3.2.* First we check the behavior of sources and targets; for the source we have

$$\begin{aligned} s\Big(\big(\mathbf{\Psi}_2\mathbf{\Psi}_1\big)(\gamma)\Big) &= s\big(\mathbf{\Psi}_2(\gamma)\big)s\big(\mathbf{\Psi}_1(\gamma)\big) = \mathbf{\Psi}_2\big(s(\gamma)\big)\mathbf{\Psi}_1\big(s(\gamma)\big) \\ &= \big(\mathbf{\Psi}_2\mathbf{\Psi}_1\big)\big(s(\gamma)\big) \end{aligned}$$

for all $\gamma \in \mathrm{Mor}(\mathbf{U})$, and similarly for targets. Next we check compositions:

$$\begin{aligned} \big(\mathbf{\Psi}_2\mathbf{\Psi}_1\big)(\gamma_2 \circ \gamma_1) &= \mathbf{\Psi}_2(\gamma_2 \circ \gamma_1)\mathbf{\Psi}_1(\gamma_2 \circ \gamma_1) \\ &= \Big(\mathbf{\Psi}_2(\gamma_2) \circ \mathbf{\Psi}_2(\gamma_1)\Big)\Big(\mathbf{\Psi}_1(\gamma_2) \circ \mathbf{\Psi}_1(\gamma_1)\Big) \\ &= \Big(\mathbf{\Psi}_2(\gamma_2)\mathbf{\Psi}_1(\gamma_2)\Big) \circ \Big(\mathbf{\Psi}_2(\gamma_1)\mathbf{\Psi}_1(\gamma_1)\Big) \end{aligned} \tag{3.27}$$

where the last equality follows on using the fundamental property of a categorical group that the group operation

$$\mathbf{G} \times \mathbf{G} \longrightarrow \mathbf{G}$$

is a functor. Hence

$$\big(\mathbf{\Psi}_2\mathbf{\Psi}_1\big)(\gamma_2 \circ \gamma_1) = \Big(\big(\mathbf{\Psi}_2\mathbf{\Psi}_1\big)(\gamma_2)\Big) \circ \Big(\big(\mathbf{\Psi}_2\mathbf{\Psi}_1\big)(\gamma_1)\Big).$$

Next, for any object $a \in \mathrm{Obj}(\mathbf{U})$ we have

$$\big(\mathbf{\Psi}_2\mathbf{\Psi}_1\big)(1_a) = \mathbf{\Psi}_2(1_a)\mathbf{\Psi}_1(1_a) = 1_{\mathbf{\Psi}_2(a)}1_{\mathbf{\Psi}_1(a)} = 1_{\mathbf{\Psi}_2(a)\mathbf{\Psi}_1(a)}.$$

Thus $\mathbf{\Psi}_2\mathbf{\Psi}_1$ is a functor.

If a functor $\mathbf{E} : \mathbf{U} \longrightarrow \mathbf{G}$ is to be the identity under pointwise multiplication its value on every object must be the identity element $e \in \mathrm{Mor}(\mathbf{G})$ and its value on every morphism must be the identity element $1_e \in \mathrm{Mor}(\mathbf{G})$. It is readily checked $\mathbf{E}$ is indeed a functor. Next, we verify that the *pointwise* inverse $\mathbf{\Psi}^{\mathrm{inv}}$ of a functor $\mathbf{\Psi}$ is also a functor. For the source we have, for any morphism $\gamma$ of $\mathbf{U}$,

$$s\Big(\mathbf{\Psi}^{\mathrm{inv}}(\gamma)\Big) = s\Big(\mathbf{\Psi}(\gamma)^{-1}\Big) = s\Big(\mathbf{\Psi}(\gamma)\Big)^{-1} = \mathbf{\Psi}^{\mathrm{inv}}\Big(s(\gamma)\Big), \tag{3.28}$$

because $s : \mathrm{Mor}(\mathbf{G}) \longrightarrow \mathrm{Obj}(\mathbf{G})$ is a homomorphism and we have denoted the multiplicative inverse of $\mathbf{\Psi}(\gamma)$ in the group $\mathrm{Mor}(\mathbf{G})$ by $\mathbf{\Psi}(\gamma)^{-1}$ (let us note that this denotes the multiplicative inverse, not compositional inverse). Thus

$$s \circ \mathbf{\Psi}^{\mathrm{inv}} = \mathbf{\Psi}^{\mathrm{inv}} \circ s.$$

Replacing $s$ by $t$ in the argument used above we see that $\mathbf{\Psi}^{\mathrm{inv}}$ also respects targets. Now consider morphisms $\gamma_1, \gamma_2 \in \mathbf{U}$ for which the composite $\gamma_2 \circ \gamma_1$ is defined. Then

$$\begin{aligned}\Big(\mathbf{\Psi}^{\mathrm{inv}}(\gamma_2) \circ \mathbf{\Psi}^{\mathrm{inv}}(\gamma_1)\Big)&\Big(\mathbf{\Psi}(\gamma_2) \circ \mathbf{\Psi}(\gamma_1)\Big)\\ &= \Big(\mathbf{\Psi}^{\mathrm{inv}}(\gamma_2)\mathbf{\Psi}(\gamma_2)\Big) \circ \Big(\mathbf{\Psi}^{\mathrm{inv}}(\gamma_1)\mathbf{\Psi}(\gamma_1)\Big)\\ &\quad \text{(since multiplication } \mathbf{G}\times\mathbf{G}\longrightarrow\mathbf{G} \text{ is a functor)}\\ &= 1_e \circ 1_e \quad \text{where } e \text{ is the identity element in } \mathrm{Obj}(\mathbf{G})\\ &= 1_e.\end{aligned} \tag{3.29}$$

Therefore,

$$\begin{aligned}\Big(\mathbf{\Psi}^{\mathrm{inv}}(\gamma_2) \circ \mathbf{\Psi}^{\mathrm{inv}}(\gamma_1)\Big) &= \Big(\mathbf{\Psi}(\gamma_2) \circ \mathbf{\Psi}(\gamma_1)\Big)^{-1}\\ &= \mathbf{\Psi}^{\mathrm{inv}}(\gamma_2 \circ \gamma_1).\end{aligned} \tag{3.30}$$

Finally,

$$\mathbf{\Psi}^{\mathrm{inv}}(1_a) = \mathbf{\Psi}(1_a)^{-1} = 1^{-1}_{\mathbf{\Psi}(a)} = 1_{\mathbf{\Psi}(a)^{-1}} = 1_{\mathbf{\Psi}^{\mathrm{inv}}(a)},$$

(with $f^{-1}$ denoting the *multiplicative* inverse of any $f \in \mathrm{Mor}(\mathbf{G})$) because

$$\mathrm{Obj}(\mathbf{G}) \longrightarrow \mathrm{Mor}(\mathbf{G}) : x \mapsto 1_x$$

is a homomorphism. Thus $\mathbf{\Psi}^{\mathrm{inv}}$ is in fact also a functor $\mathbf{U} \longrightarrow \mathbf{G}$.

Finally, it is clear that pointwise multiplication on functors $\mathbf{U} \longrightarrow \mathbf{G}$ is an associative operation. □

*Proof of Proposition 3.3.* For any $a \in \mathrm{Obj}(\mathbf{U})$, the values $\mathbf{\Psi}_1(a)$ and $\mathbf{\Psi}_2(a)$ are objects of $\mathbf{G}$, which means that they are elements of $G$, and $\mathbb{T}(a) : \mathbf{\Psi}_1(a) \longrightarrow \mathbf{\Psi}_2(a)$, being a morphism of $\mathbf{G}$, is given by some $\big(h_{\mathbb{T}}(a), g_{\mathbb{T}}(a)\big) \in H \rtimes_\alpha G = \mathrm{Mor}(\mathbf{G})$:

$$\mathbb{T}(a) = \big(h_{\mathbb{T}}(a), g_{\mathbb{T}}(a)\big) \in H \rtimes_\alpha G. \tag{3.31}$$

Then

$$\begin{aligned}\mathbf{\Psi}_1(a) &= s\big(h_{\mathbb{T}}(a), g_{\mathbb{T}}(a)\big) = g_{\mathbb{T}}(a)\\ \mathbf{\Psi}_2(a) &= t\big(h_{\mathbb{T}}(a), g_{\mathbb{T}}(a)\big) = \tau\big(h_{\mathbb{T}}(a)\big)g_{\mathbb{T}}(a).\end{aligned} \tag{3.32}$$

This proves (3.11).

Next, working with any morphism

$$f : a \longrightarrow b$$

in $\mathbf{U}$, we write

$$\mathbf{\Psi}_i(f) = \big(h_i(f), g_i(f)\big) \qquad \text{for } i \in \{1, 2\},$$

where the source $g_i(f)$ is, as we have already observed in (3.32), given by $g_{\mathbb{T}}(a)$ if $i = 1$ and by $\tau\big(h_{\mathbb{T}}(a)\big)g_{\mathbb{T}}(a)$ if $i = 2$. Using the commuting diagram (3.10), expressed more explicitly as

$$\mathbf{\Psi}_2(f) \circ \mathbb{T}(a) = \mathbb{T}(b) \circ \mathbf{\Psi}_1(f), \tag{3.33}$$

we have

$$\big(h_2(f), g_2(f)\big) \circ \mathbb{T}(a) = \mathbb{T}(b) \circ \mathbf{\Psi}_1(f). \tag{3.34}$$

Using the expression $\mathbb{T}(a) = \big(h_{\mathbb{T}}(a), g_{\mathbb{T}}(a)\big)$ and similarly for $\mathbb{T}(b)$ we then have

$$\big(h_2(f)h_{\mathbb{T}}(a), g_{\mathbb{T}}(a)\big) = \big(h_{\mathbb{T}}(b)h_1(f), g_{\mathbb{T}}(a)\big). \tag{3.35}$$

From this we see that

$$h_2(f) = h_{\mathbb{T}}(b)h_1(f)h_{\mathbb{T}}(a)^{-1},$$

which is the same as (3.12). $\square$

## 4. Sections and Product Structure

As before, we work with a principal categorical bundle

$$\pi : \mathbf{P} \longrightarrow \mathbf{B},$$

having structure categorical group $\mathbf{G}$. Let $\mathbf{U}$ be a subcategory of $\mathbf{B}$. A *section* $\sigma$ of the bundle $\pi : \mathbf{P} \longrightarrow \mathbf{B}$ *over* $\mathbf{U}$ is a functor

$$\sigma : \mathbf{U} \longrightarrow \mathbf{P}$$

for which $\pi \circ \sigma$ is the identity functor. Now let

$$\mathbf{P}_{\mathbf{U}} \stackrel{\text{def}}{=} \pi^{-1}(\mathbf{U})$$

be the subcategory of $\mathbf{P}$ whose objects and morphisms project down to objects and morphisms, respectively, of $\mathbf{U}$. Then the functor

$$\pi_{\mathbf{U}} : \mathbf{P}_{\mathbf{U}} \longrightarrow \mathbf{U}$$

given by restriction of $\pi$ is a principal categorical bundle with structure categorical group $\mathbf{G}$, as can be readily verified. Consider then

$$\mathbf{\Psi}_\sigma : \mathbf{U} \times \mathbf{G} \longrightarrow \mathbf{P}_{\mathbf{U}} \tag{4.1}$$

defined on objects by

$$\mathbf{\Psi}_\sigma(a, g) = \sigma(a)g \qquad \text{for all } a \in \mathrm{Obj}(\mathbf{U}) \text{ and } g \in \mathrm{Obj}(\mathbf{G}), \tag{4.2}$$

and on morphisms by

$$\mathbf{\Psi}_\sigma\big(\gamma, \psi\big) = \sigma(\gamma)\psi \qquad \text{for all } \gamma \in \mathrm{Mor}(\mathbf{U}) \text{ and } \psi \in \mathrm{Mor}(\mathbf{G}). \tag{4.3}$$

It is apparent that $\boldsymbol{\Psi}_\sigma$ is equivariant with respect to the action of $\mathbf{G}$ both on objects and on morphisms; for example, at the object level,

$$\begin{aligned} \boldsymbol{\Psi}_\sigma\big((\gamma,\psi)\psi_1\big) &= \boldsymbol{\Psi}_\sigma(\gamma,\psi\psi_1) \\ &= \sigma(\gamma)\psi\psi_1 \\ &= \boldsymbol{\Psi}_\sigma(\gamma,\psi)\psi_1. \end{aligned} \tag{4.4}$$

Moreover,

$$\pi_{\mathbf{U}} \circ \boldsymbol{\Psi}_\sigma = p_{\mathbf{U}} : \mathbf{U} \times \mathbf{G} \longrightarrow \mathbf{U},$$

which is the projection on the first factor. The conditions (b2) and (b3) defining a categorical principal bundle require that the action of $\mathbf{G}$ be free (conditon (b2)) and transitive (conditon (b3)) on each fiber both for objects and for morphisms. These, along with the observations we have just made, imply that $\boldsymbol{\Psi}_\sigma$ is bijective, again both on objects and on morphisms. Let us state the argument for objects (the case for morphisms being exactly analogous). First, if $p \in P$ and $a = \pi(p)$ then $p$ and $\sigma(a)$ lie on the same fiber and so $p = \sigma(a)g$ for some $g \in G$; thus $\boldsymbol{\Psi}_\sigma$ is surjective on objects. Next, suppose $\boldsymbol{\Psi}_\sigma(a,g) = \boldsymbol{\Psi}_\sigma(a',g')$; applying $\pi$ we have

$$a = \pi\big(\boldsymbol{\Psi}_\sigma(a,g)\big) = \pi\big(\boldsymbol{\Psi}_\sigma(a',g')\big) = a',$$

and then from

$$\boldsymbol{\Psi}_\sigma(a)g = \boldsymbol{\Psi}_\sigma(a)g'$$

it follows that $g = g'$ because the action of $G$ is free. Thus $\boldsymbol{\Psi}_\sigma$ is also injective.

Finally, let us verify that $\boldsymbol{\Psi}_\sigma$ preserves composition of morphisms. Consider the composite morphism

$$(\gamma_2,\psi_2) \circ (\gamma_1,\psi_1) \in \mathrm{Mor}(\mathbf{U} \times \mathbf{G});$$

then

$$\begin{aligned} \boldsymbol{\Psi}_\sigma\Big((\gamma_2,\psi_2)\circ(\gamma_1,\psi_1)\Big) &= \boldsymbol{\Psi}_\sigma(\gamma_2\circ\gamma_1, \psi_2\circ\psi_1) \\ &= \sigma(\gamma_2\circ\gamma_1)(\psi_2\circ\psi_1) \\ &= \Big(\sigma(\gamma_2)\circ\sigma(\gamma_1)\Big)(\psi_2\circ\psi_1) \\ &\qquad \text{(because } \sigma \text{ is a functor)} \\ &= \Big(\sigma(\gamma_2)\psi_2\Big)\circ\Big(\sigma(\gamma_1)\psi_1\Big) \\ &\text{(because the action } \mathbf{P}\times\mathbf{G}\longrightarrow\mathbf{P} \text{ is a functor)} \\ &= \boldsymbol{\Psi}_\sigma(\gamma_2,\psi_2)\circ\boldsymbol{\Psi}_\sigma(\gamma_1,\psi_1). \end{aligned} \tag{4.5}$$

To summarize, we have proved the ‘if’ part of the categorical analog of a basic observations about bundles:

**Proposition 4.1.** *Suppose $\pi : \mathbf{P} \longrightarrow \mathbf{U}$ is a principal categorical bundle with structure categorical group $\mathbf{G}$. Then $\pi$ admits a section if and only if there is an isomorphism of categories*

$$\boldsymbol{\Psi} : \mathbf{U} \times \mathbf{G} \longrightarrow \mathbf{P}$$

*that is equivariant with respect to the action of $\mathbf{G}$ and preserves fibers in the sense that*

$$p_U \circ \boldsymbol{\Psi} = \pi,$$

*where $p_U : \mathbf{U} \times \mathbf{G} \longrightarrow \mathbf{U}$ is the functor that projects on the first factor.*

The 'only if' part is readily verified: it is clear that

$$\mathbf{U} \longrightarrow \mathbf{U} \times \mathbf{G} : u \mapsto (u, e),$$

with obvious meaning, is a section of $p_{\mathbf{U}}$.

Let us now consider the case in the preceding Proposition when $\mathbf{P}$ is itself also the product categorical bundle $\mathbf{U} \times \mathbf{G}$; thus consider a categorical principal bundle isomorphism functor

$$\Phi : \mathbf{U} \times \mathbf{G} \longrightarrow \mathbf{U} \times \mathbf{G}; \tag{4.6}$$

it is equivariant with respect to the $\mathbf{G}$-action and preserves the first component. Consequently, it is given on objects and morphisms by

$$\begin{aligned} \Phi(a, g) &= \Phi(a, e)g = \big(a, \sigma_\Phi(a)g\big) \\ &\quad \text{for all } (a, g) \in \mathrm{Obj}(U) \times \mathrm{Obj}(\mathbf{G}) \\ \Phi(\gamma, \psi) &= \Phi(\gamma, e)\phi = \big(\gamma, \sigma_\Phi(\gamma)\phi\big), \\ &\quad \text{for all } (\gamma, \phi) \in \mathrm{Mor}(U) \times \mathrm{Mor}(\mathbf{G}), \end{aligned} \tag{4.7}$$

where $e$ in the first line is the identity in $\mathrm{Obj}(\mathbf{G})$ and in the second line it is the identity in $\mathrm{Mor}(\mathbf{G})$. Functoriality of $\Phi$ means, in particular, that it preserves sources and targets in the sense that $\Phi s = s\Phi$ and $\Phi t = t\Phi$; applying this in (4.7) with $g = e$ and $\psi = e = 1_e$ we have

$$\begin{aligned} s\sigma_\Phi(\gamma) &= \sigma_\Phi s(\gamma) \\ t\sigma_\Phi(\gamma) &= \sigma_\Phi t(\gamma), \end{aligned} \tag{4.8}$$

where $\sigma_\Phi$ is as given by (4.7). It is also clear from the functoriality of $\Psi$ that $\sigma_\Psi$ takes compositions to compositions and identity morphisms to identity morphisms. Thus $\sigma_\Phi$ is a functor $\mathbf{U} \longrightarrow \mathbf{G}$. We summarize this observation now.

**Proposition 4.2.** *Suppose $\mathbf{U}$ is a category, $\mathbf{G}$ a categorical group, and*

$$\Phi : \mathbf{U} \times \mathbf{G} \longrightarrow \mathbf{U} \times \mathbf{G}; \tag{4.9}$$

*is a functor that preserves the first component and is equivariant under the right action of* $\mathbf{G}$ *on* $\mathbf{U} \times \mathbf{G}$*. Then there is a functor*

$$\sigma_\Phi : \mathbf{U} \longrightarrow \mathbf{G} \tag{4.10}$$

*such that* $\Phi$ *is given on objects and morphisms by (4.7). Conversely, if* $\sigma : \mathbf{U} \longrightarrow \mathbf{G}$ *is a functor then* $\Psi$*, as specified by (4.7) with* $\sigma$ *replacing* $\sigma_\Psi$*, is a functor that preserves the first component and is* $\mathbf{G}$*-equivariant.*

*Proof.* We have already established the direct implication. The converse is readily verified. □

Continuing with the notation $\sigma_\Psi$ we can see that if

$$\boldsymbol{\Psi}_1, \boldsymbol{\Psi}_2 : \mathbf{U} \times \mathbf{G} \longrightarrow \mathbf{U} \times \mathbf{G}$$

are functors then $\boldsymbol{\Psi}_2 \circ \boldsymbol{\Psi}_1$ is also a functor. The corresponding functor $\mathbf{U} \longrightarrow \mathbf{G}$ is given by

$$\sigma_{\boldsymbol{\Psi}_2 \circ \boldsymbol{\Psi}_1} = \sigma_{\boldsymbol{\Psi}_2} \sigma_{\boldsymbol{\Psi}_1}. \tag{4.11}$$

Thus the product operation in the categorical group $\mathbf{G}^{\mathbf{U}}$ corresponds to composition of functorial automorphisms of $\mathbf{U} \times \mathbf{G}$.

## 5. Functorial cocycles

We continue working with a principal categorical bundle

$$\pi : \mathbf{P} \longrightarrow \mathbf{B}$$

with structure categorical group $\mathbf{G}$.

Now suppose that the 'object' bundle

$$\pi : P = \mathrm{Obj}(\mathbf{P}) \longrightarrow B = \mathrm{Obj}(\mathbf{B})$$

is a principal $G$-bundle in the usual sense. In particular, $P$ and $B$ are topological spaces.

Let

$$\{U_i\}_{i \in I}$$

be an open covering of $B$ for which the bundle $\pi : P \longrightarrow B$ trivializes over each $U_i$ in the sense that there is a $G$-equivariant diffeomorphism

$$\phi_i : U_i \times G \longrightarrow \pi^{-1}(U_i) \tag{5.1}$$

for which any point of the form $(u, g)$ is mapped to a point on the fiber in $P$ lying over $u \in U_i$.

We denote by

$$\mathbf{U}_i^j$$

the category whose objects are the points of $U_i$ labeled with $i$ along with the points of $U_j$, labeled with $j$, and whose morphisms are the morphisms $\gamma$ of $\mathbf{B}$ with source in $U_i$ and target in $U_j$, along with all the needed identity morphisms. Thus

$$\begin{aligned}
&\mathrm{Obj}(\mathbf{U}_i^j) = (\{i\} \times U_i) \cup (\{j\} \times U_j) \\
&\mathrm{Mor}(\mathbf{U}_i^j) \\
&\quad = \{\gamma \in \mathrm{Mor}(\mathbf{B}) : s(\gamma) \in U_i, t(\gamma) \in U_j\} \cup \{\mathrm{id}_x \,:\, x \in \mathrm{Obj}(\mathbf{U}_i^j)\}.
\end{aligned} \tag{5.2}$$

Source and target mappings are defined by requiring that that source of any $\gamma \in \mathrm{Mor}(\mathbf{U}_i^j)$ be $\big(j, s(\gamma)\big)$ and the target be $\big(j, t(\gamma)\big)$. Note that, technically, these source and target maps depend on $i$ and $j$. In the object set $\mathrm{Obj}(\mathbf{U}_i^j)$ we are labeling the points of $U_i$ with $i$, a tracking mechanism that is useful when considering transition functions. Figure 1 illustrates a typical morphism.

In the notation used for the categories $\mathbf{U}_i^k$ *subscripts indicate the source* and *superscripts indicate the target*. We will maintain this convention in our notation for categories whose morphisms begin and end in potentially different subsets of a given object space.

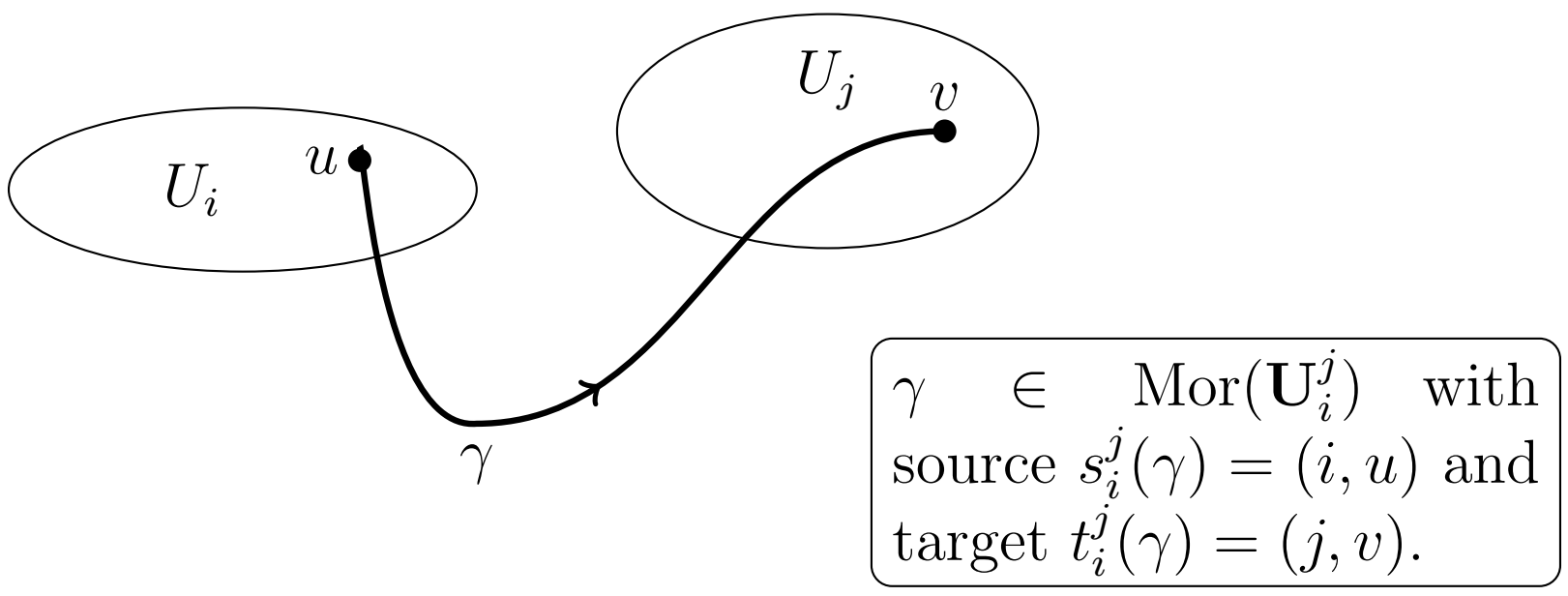

FIGURE 1. A morphism of $\mathbf{U}_i^j$

If we elevate $\phi_i$ to a local product strucure at the categorical level then we should have, by analogy, an isomorphism of categorical principal bundles

$$\mathbf{\Phi}_i^j : \mathbf{U}_i^j \times \mathbf{G} \longrightarrow \mathbf{P}_i^j, \tag{5.3}$$

where $\mathbf{P}_i^j$ is the part of $\mathbf{P}$ that projects down to $\mathbf{U}_i^j$. Let us be more precise about what this means: the objects of $\mathbf{P}_i^j$ are of the form $(i, p)$ with $p \in \pi^{-1}(U_i)$ and of the form $(j, q)$ with $q \in \pi^{-1}(U_j)$; a morphism of $\mathbf{P}_i^j$ is either an identity morphism or a morphism $\gamma \in \mathrm{Mor}(\mathbf{P})$ whose source lies in $\pi^{-1}(U_i)$ and whose target lies in $\pi^{-1}(U_j)$. Thus

$$\mathrm{Obj}(\mathbf{P}_i^j) = \{(i, p) : p \in \pi^{-1}(U_i)\} \cup \{(j, q) : q \in \pi^{-1}(U_j)\}, \tag{5.4}$$

and morphisms are given similarly. The source and target maps are defined in the obvious way.

Let us now consider what happens on the overlap category

$$\mathbf{U}_{ik}^{jl}, \tag{5.5}$$

whose object set and morphism set are given as follows. The object set of $\mathbf{U}_{ik}^{jl}$ consists of all points of the form $(ik, u)$ with $u \in U_i \cap U_k$, along with all points of the form $(jl, v)$ with $v \in U_j \cap U_l$:

$$\mathrm{Obj}(\mathbf{U}_{ik}^{jl}) = \{(ik, u) \,:\, u \in U_i \cap U_k\} \cup \{(jl, v) \,:\, v \in U_j \cap U_l\}. \tag{5.6}$$

Aside from the identity morphisms, a morphism $\gamma$ of $\mathbf{U}_{ik}^{jl}$ is a morphism of $\mathbf{B}$ with source in $U_i \cap U_k$ and target in $U_j \cap U_l$; the source of $\gamma$ in $\mathbf{U}_{ik}^{jl}$ is then $\big(ik, s(\gamma)\big)$ and the target is $\big(jl, t(\gamma)\big)$. (More technically, the 'intersection' (5.5) is a product in the category given by the partially ordered set whose elements are of the form $\mathbf{U}_{\alpha}^{\beta}$, with $\alpha$ and $\beta$ being $n$-tuples of indices and $n \geq 1$, and the ordering relation is obtained by inclusions.) Figure 2 illustrates the idea.

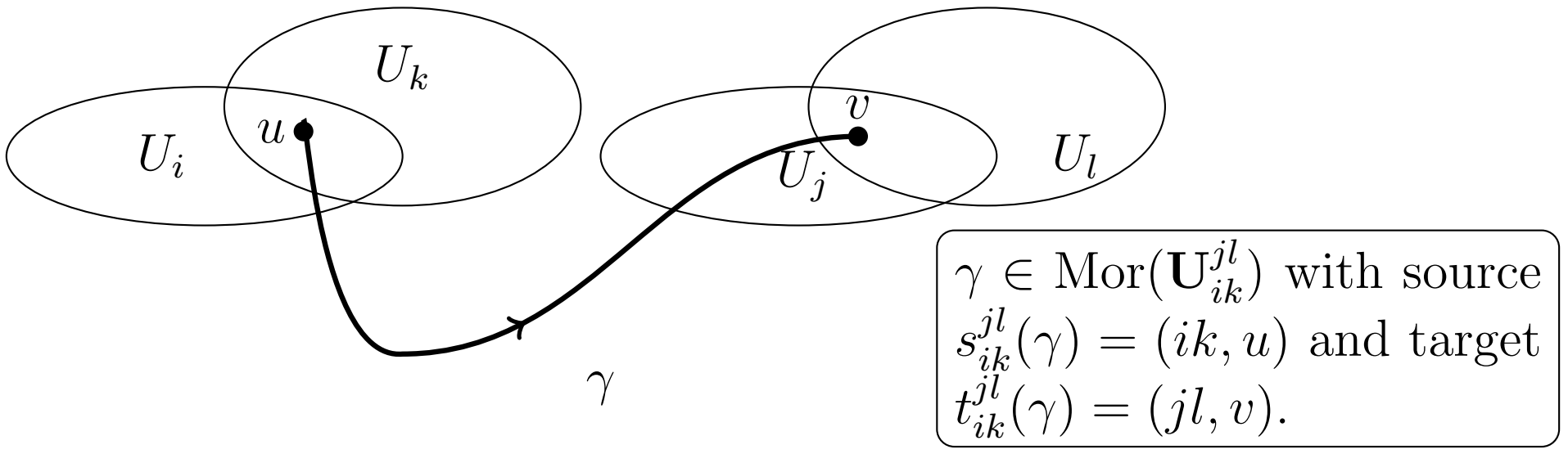

FIGURE 2. A morphism of $\mathbf{U}_{ik}^{jl}$

A local product structure of $\mathbf{P}_k^l$ is expressed through a functorial isomorphism of categorical bundles

$$\mathbf{\Phi}_k^l : \mathbf{U}_k^l \times \mathbf{G} \longrightarrow \mathbf{P}_k^l. \tag{5.7}$$

Transitioning from this trivialization of $\mathbf{P}_k^l$ to the corresponding trivialization of $\mathbf{P}_i^j$, assuming that $U_i \cap U_k$ and $U_j \cap U_l$ are nonempty, we obtain the functor

$$\mathbf{\Sigma}_{ik}^{jl} : \mathbf{U}_{ik}^{jl} \times \mathbf{G} \longrightarrow \mathbf{U}_{ik}^{jl} \times \mathbf{G}, \tag{5.8}$$

obtained by first applying $\mathbf{\Phi}_k^l$ and then applying the inverse of $\mathbf{\Phi}_i^j$. Note that the functor $\mathbf{\Sigma}_{ik}^{jl}$ describes transition *from* the trivialization labeled $\Big\{\Big\}_k^l$ *to* the one labeled $\Big\{\Big\}_i^j$.

By the method of Proposition 4.2 (specifically (4.7)) this corresponds then to a functor

$$\sigma^{j\,l}_{i\,k} : \mathbf{U}^{jl}_{ik} \longrightarrow \mathbf{G}. \tag{5.9}$$

For example, at the object level this is obtained from the functor $\mathbf{\Sigma}^{j\,l}_{i\,k}$ by picking out the second component:

$$\mathbf{\Sigma}^{j\,l}_{i\,k}\big(x, e\big) = \big(x, \sigma^{j\,l}_{i\,k}(x)\big) \tag{5.10}$$

where $x$ is either of the form $(ik, u)$ with $u \in U_i \cap U_k$ or of the form $(jl, v)$ with $v \in U_j \cap U_l$. Taking $x$ to be a morphism $\gamma$ of $\mathbf{U}^{jl}_{ik}$ and $e$ the identity in $\mathrm{Mor}(\mathbf{G})$ produces the definition of $\sigma^{j\,l}_{i\,k}(\gamma)$.

Thus,

$$\mathbf{\Phi}^{j\,-1}_{i}\Big(\mathbf{\Phi}^{l}_{k}\big(x, e\big)\Big) = (x, e)\sigma^{j\,l}_{i\,k}(x) \tag{5.11}$$

for all $x$ for which the expression on the left makes sense. This is the exact functorial counterpart of the definition of the classical transition function $g_{\alpha\beta}$ given in (1.2).

Even more explicitly, taking $x$ to be an object of the form $(ik, u)$ we first apply $\mathbf{\Phi}^{l}_{k}$ to $\big((k, u), e\big)$ to obtain an object

$$(k, p^*) = \mathbf{\Phi}^{l}_{k}\big((k, u), e\big) \in \mathrm{Obj}(\mathbf{P}^{l}_{k}) \tag{5.12}$$

and then take $\sigma^{j\,l}_{i\,k}(x)$ to be the element of the group $\mathrm{Obj}(\mathbf{G})$ for which

$$\mathbf{\Phi}^{j}_{i}\big((i, u), \sigma^{j\,l}_{i\,k}(x)\big) = (k, p^*). \tag{5.13}$$

Thus at the object level the transition functor $\sigma^{j\,l}_{i\,k}$ is given by an $\mathrm{Obj}(\mathbf{G})$-valued function on $\mathrm{Obj}(\mathbf{U}^{k}_{j})$:

$$\{(ik, u) : u \in U_i \cap U_k\} \cup \{(jl, v) : v \in U_j \cap U_l\} \longrightarrow \mathrm{Obj}(\mathbf{G}), \tag{5.14}$$

a counterpart to the traditional function that is made up of a bundle transition function over $U_i \cap U_k$ along with a bundle transition function over $U_j \cap U_l$.

We turn now to discussion the behavior of transition functors in a triple overlap of the categories $\mathbf{U}^{j}_{i}$. For this we need to first specify some necessary notation for triple overlaps.

We define triple overlap categories $\mathbf{U}^{jln}_{ikm}$ entirely analogously to the double overlaps we have considered before. We take the object set of $\mathbf{U}^{jln}_{ikm}$ to be $U_{ikm} \cup U_{jln}$, where a typical point of $U_{ikm}$ is of the form $(ikm, a)$ with $a \in U_i \cap U_k \cap U_m$. Morphisms are all the identity morphisms along with the morphisms of $\mathbf{U}$ that have source in $U_i \cap U_k \cap U_m$ and target in $U_j \cap U_l \cap U_n$, as illustrated in Figure 3.

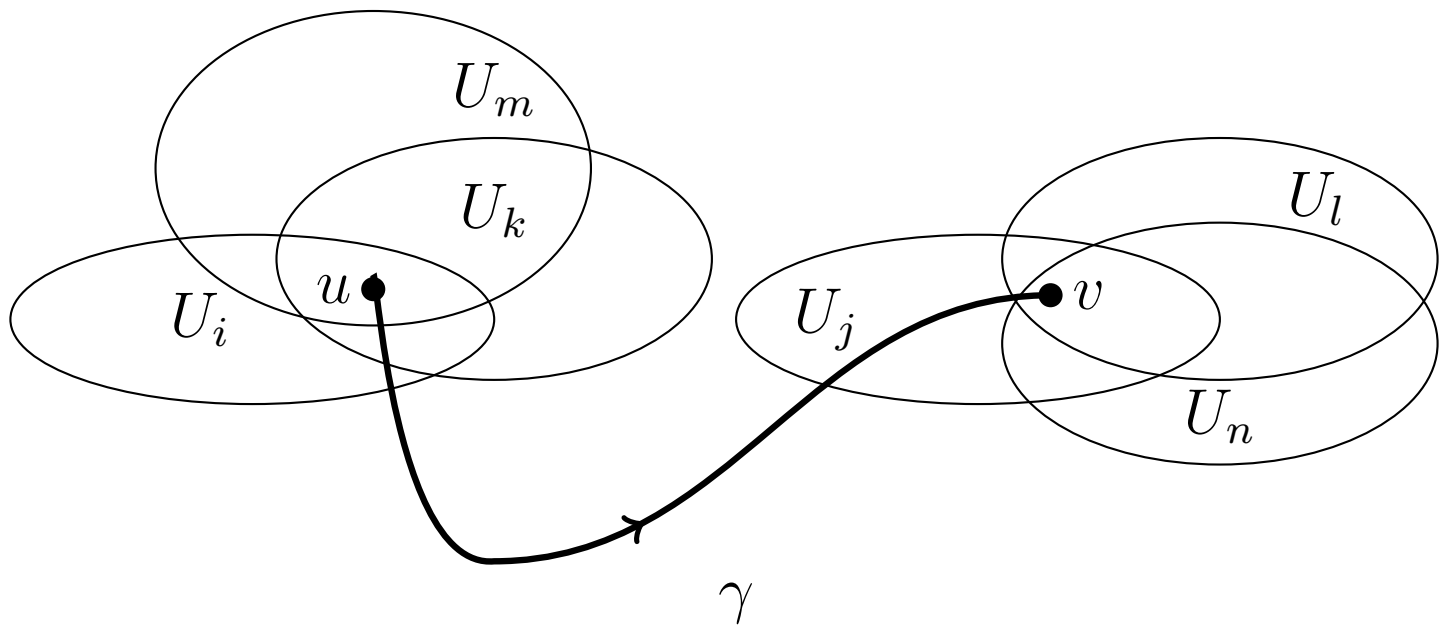

FIGURE 3. A morphism $\gamma \in \mathrm{Mor}(\mathbf{U}^{jln}_{ikm})$ with source $(ikm, u)$ and target $(jln, v)$.

Restricting the transition functor $\sigma^{j\,l}_{i\,k}$ to the triple overlap category $\mathbf{U}^{jln}_{ikm}$ we obtain

$$\sigma^{j\,l\,|n|}_{i\,k\,|m|} : \mathbf{U}^{jln}_{ikm} \longrightarrow \mathbf{G}. \tag{5.15}$$

Thus for an object $(ikm, u)$ of $\mathbf{U}^{jln}_{ikm}$ we set

$$\begin{aligned} \sigma^{j\,l\,|n|}_{i\,k\,|m|}(ikm, u) &= \sigma^{j\,l}_{i\,k}(ik, u) \\ \sigma^{j\,l\,|n|}_{i\,k\,|m|}(jln, v) &= \sigma^{j\,l}_{i\,k}(jl, v), \end{aligned} \tag{5.16}$$

and for a morphism $\gamma \in \mathrm{Mor}(\mathbf{U}^{jln}_{ikm})$ we set

$$\sigma^{j\,l\,|n|}_{i\,k\,|m|}(\gamma) = \sigma^{j\,l}_{i\,k}(\gamma). \tag{5.17}$$

Restricting the other transition functors $\sigma^{\cdot\cdot}_{\cdot\cdot}$ in different ways we have the functor

$$\sigma^{|j|\,l\,n}_{|i|\,k\,m} : \mathbf{U}^{jln}_{ikm} \longrightarrow \mathbf{G} \tag{5.18}$$

and the functor

$$\sigma^{j\,|l|\,n}_{i\,|k|\,m} : \mathbf{U}^{jln}_{ikm} \longrightarrow \mathbf{G}. \tag{5.19}$$

Consistency of transitions requires that

$$\sigma^{j\,l\,|n|}_{i\,k\,|m|}\sigma^{|j|\,l\,n}_{|i|\,k\,m} = \sigma^{j\,|l|\,n}_{i\,|k|\,m}. \tag{5.20}$$

This equation is more easily understood by simplifying the notation:

$$\sigma^{j\,l}_{i\,k}\sigma^{l\,n}_{k\,m} = \sigma^{j\,n}_{i\,m}, \tag{5.21}$$

where it is understood that each functor $\sigma^{\cdot\cdot}_{\cdot\cdot}$ has been restricted to the triple overlap $\mathbf{U}^{jln}_{ikm}$. This is just as in the classical cocycle relation (1.3):

$$g_{\alpha\beta}g_{\beta\gamma} = g_{\alpha\gamma} \tag{5.22}$$

for transition functions $g_{\alpha\beta}$ defined on $U_\alpha \cap U_\beta$ , and in this relation it is understood that the functions $g_{\cdot\cdot}$ are all restricted to $U_\alpha \cap U_\beta \cap U_\gamma$. The analogy between the classical cocycle relation (5.22) and the functorial one (5.21) is seen more sharply on thinking of the index $\alpha$ as $\binom{j}{i}$, the index $\beta$ as $\binom{l}{k}$ and $\gamma$ as $\binom{n}{m}$.

To establish the cocycle relation (5.20) we apply the relation

$$\big((\mathbf{\Phi}^j_i)^{-1}\circ \mathbf{\Phi}^l_k\big)\circ\big((\mathbf{\Phi}^l_k)^{-1}\circ \mathbf{\Phi}^n_m\big) = (\mathbf{\Phi}^j_i)^{-1}\circ \mathbf{\Phi}^n_m \tag{5.23}$$

to $(x,e)$ which is either an object or a morphism of $\mathbf{U}^{jln}_{ikm}$; however, in order to apply the double-indexed functors $\sigma^{\cdot\cdot}_{\cdot\cdot}$ to $x$ we need to restrict them first, for example as

$$\sigma^{|j|\,l\,n}_{|i|\,k\,m}(x).$$

Thus the right side of (5.23) applied to $(x,e)$ gives:

$$(\mathbf{\Phi}^j_i)^{-1}\circ \mathbf{\Phi}^n_m(x,e) = \big(x, \sigma^{j\,|l|\,n}_{i\,|k|\,m}(x)\big). \tag{5.24}$$

We apply the left side of (5.23) to $(x,e)$ by repeating this procedure:

$$\begin{aligned}\big((\mathbf{\Phi}^j_i)^{-1}\circ \mathbf{\Phi}^l_k\big)\circ\big((\mathbf{\Phi}^l_k)^{-1}\circ \mathbf{\Phi}^n_m\big)(x,e) &= \big((\mathbf{\Phi}^j_i)^{-1}\circ \mathbf{\Phi}^l_k\big)\big(x,\sigma^{|j|\,l\,n}_{|i|\,k\,m}(x)\big)\\ &= \big((\mathbf{\Phi}^j_i)^{-1}\circ \mathbf{\Phi}^l_k\big)\big(x,e\big)\sigma^{|j|\,l\,n}_{|i|\,k\,m}(x)\\ &= \big(x,\sigma^{j\,l}_{i\,k}(x)\big)\sigma^{|j|\,l\,n}_{|i|\,k\,m}(x)\\ &= (x,e)\sigma^{j\,l\,|n|}_{i\,k\,|m|}(x)\sigma^{|j|\,l\,n}_{|i|\,k\,m}(x).\end{aligned} \tag{5.25}$$

This establishes the consistency relation (5.20) among transition functors.

The system of functors $\{\sigma^{\cdot\cdot}_{\cdot\cdot}\}$ forms a special type of *functorial cocycle.* We will return to this shortly, after exploring $\mathbf{G}$-valued cocycles directly, without involving bundles in the definition.

5.1. **Cocycles.** Consider as given the following data:

(gh2) for every $i,j\in I$ a smooth function

$$h_{ij}: U_i\cap U_j \longrightarrow H$$

(gh3) for every $i,j,k\in I$ a smooth function

$$h_{ijk}: U_i\cap U_j\cap U_k \longrightarrow H$$

satisfying the following cocycle condition

$$h_{ijk}h_{ik} = h_{ij}h_{jk} \qquad \text{on } U_i\cap U_j\cap U_k, \tag{5.26}$$

for all $i,j,k\in I$; thus applying $h_{ijk}$ to $h_{ik}$ has the effect of splitting it into two factors with the subscripts $ij$ and $jk$.

As we have seen above, it will sometimes be necessary to label a point of $U_i \cap U_j$ with $ij$ to encode the information about the indices $i$ and $j$; to this end we write

$$U_{ij} = \{ij\} \times (U_i \cap U_j), \tag{5.27}$$

and define $U_{ijk}$ analogously. Moreover, we use the notation

$$U_i^j = \{\gamma \in \mathrm{Mor}(\mathbf{B}) : s(\gamma) \in U_i, t(\gamma) \in U_j\}. \tag{5.28}$$

and

$$U_{ik}^{jl} \qquad \text{and} \qquad U_{ikm}^{jln},$$

defined analogously; for example,

$$U_{ik}^{jl} = \{\gamma \in \mathrm{Mor}(\mathbf{B}) : s(\gamma) \in U_i \cap U_k \quad \text{and} \qquad t(\gamma) \in U_j \cap U_l\}. \tag{5.29}$$

Let us recall from (5.5) the categories $\mathbf{U}_{ik}^{jl}$ specified by:

$$\mathrm{Obj}(\mathbf{U}_{ik}^{jl}) = U_{ik} \cup U_{jl} \tag{5.30}$$

and

$$\mathrm{Mor}(\mathbf{U}_{ik}^{jl}) = U_{ik}^{jl} \cup \{\mathrm{id}_a \,:\, a \in U_{jl} \cup U_{ik}\}, \tag{5.31}$$

where $\mathrm{id}_a$ is the identity morphism $a \longrightarrow a$. The source and target maps are defined in the obvious way: if $\gamma \in \mathrm{Mor}(\mathbf{U}_{ik}^{jl})$ then

$$s_{ik}^{jl}(\gamma) = \big(ik, s(\gamma)\big) \qquad \text{and} \qquad t_{ik}^{jl}(\gamma) = \big(jl, t(\gamma)\big). \tag{5.32}$$

Thus the overlap category $\mathbf{U}_{ik}^{jl}$ is 'almost' $\mathbf{U}_i^j \cap \mathbf{U}_k^l$, except that the objects are labeled with $ik$ or $jl$; or else, if for example $i \neq k$, the intersection $\mathbf{U}_i^j \cap \mathbf{U}_k^l$ would be empty.

5.2. **Functorial cocycles.** Now let us define functors

$$\theta_{i\,k}^{j\,l} : \mathbf{U}_{ik}^{jl} \longrightarrow \mathbf{G} \tag{5.33}$$

given on objects by

$$\theta_{i\,k}^{j\,l}(mn, a) = g_{mn}(a) \stackrel{\text{def}}{=} \tau\big(h_{mn}(a)\big) \tag{5.34}$$

for all $a \in U_{mn}$, where $mn \in \{ik, jl\}$, and on morphisms by

$$\theta_{i\,k}^{j\,l}(\gamma) = \Big(h_{jl}(\gamma_1) h_{ik}(\gamma_0)^{-1}, g_{ik}(\gamma_0)\Big) \tag{5.35}$$

where $\gamma_0 = s(\gamma)$ and $\gamma_1 = t(\gamma)$. When $\gamma = \mathrm{id}_a$ we take

$$\theta_{i\,k}^{j\,l}(\mathrm{id}_a) = \Big(e, g_{ik}(a)\Big) \tag{5.36}$$

The category $\mathbf{U}_{ikm}^{jln}$ can be viewed as a subcategory of $\mathbf{U}_{im}^{jn}$ in the obvious way (with objects identified using the appropriate indices).

On the triple overlap category $\mathbf{U}_{ikm}^{jln}$, assumed non-empty, we have the "restricted" functors (as in (5.16))

$$\begin{aligned}\theta_{i\,k\,|m|}^{j\,l\,|n|} &= \theta_{ik}^{jl}|\mathbf{U}_{ikm}^{jln} : \mathbf{U}_{ikm}^{jln} \longrightarrow \mathbf{G}\\ \theta_{|i|\,k\,m}^{|j|\,l\,n} &= \theta_{km}^{ln}|\mathbf{U}_{ikm}^{jln} : \mathbf{U}_{ikm}^{jln} \longrightarrow \mathbf{G}.\end{aligned} \tag{5.37}$$

We also have the restricted functor

$$\theta_{i\,|k|\,m}^{j\,|l|\,n} = \theta_{im}^{jn}|\mathbf{U}_{ikm}^{jln} : \mathbf{U}_{ikm}^{jln} \longrightarrow \mathbf{G}. \tag{5.38}$$

There is a fundamental consistency relation between these restrictions, as we explain in the following result.

**Proposition 5.1.** *With framework and notation as above, there is a natural transformation*

$$\mathbb{T} : \big(\theta_{i\,k\,|m|}^{j\,l\,|n|}\big)\big(\theta_{|i|\,k\,m}^{|j|\,l\,n}\big) \longrightarrow \theta_{i\,|k|\,m}^{j\,|l|\,n} \tag{5.39}$$

*for which the function $h_{\mathbb{T}}$ of Proposition 3.3 is given on the overlap object set $U_{ikm}$ by $(ikm, a) \mapsto h_{ikm}(a)$ and on $U_{jln}$ by $(jln, a) \mapsto h_{jln}(a)$.*

Comparing with (5.20), we see that in the earlier context $\mathbb{T}$ is simply equality in that case. We can also state (5.39) as we did (5.21) by simplifying the notation:

$$\mathbb{T} : \theta_{i\,k}^{j\,l}\theta_{k\,m}^{l\,n} \longrightarrow \theta_{i\,m}^{j\,n}. \tag{5.40}$$

*Proof.* We work with any $\gamma \in U_{ikm}^{jln}$. We have

$$\begin{aligned}\theta_{i\,k\,|m|}^{j\,l\,|n|}(\gamma) &= \Big(h_{jl}(\gamma_1)h_{ik}(\gamma_0)^{-1}, g_{ik}(\gamma_0)\Big)\\ \theta_{|i|\,k\,m}^{|j|\,l\,n}(\gamma) &= \Big(h_{ln}(\gamma_1)h_{km}(\gamma_0)^{-1}, g_{km}(\gamma_0)\Big)\end{aligned} \tag{5.41}$$

Using the formula (2.12) we see that the $H$-component of the product $\big(\theta_{i\,k\,|m|}^{j\,l\,|n|}\big)\big(\theta_{|i|\,k\,m}^{|j|\,l\,n}\big)$ is

$$h_{jl}(\gamma_1)h_{ik}(\gamma_0)^{-1}g_{ik}(\gamma_0)\ \ h_{ln}(\gamma_1)h_{km}(\gamma_0)^{-1}g_{ik}(\gamma_0)^{-1}. \tag{5.42}$$

'Gauge-transforming' by the values of the function $h_{ikm}$ and $h_{jln}$ at the source and target points, respectively, produces:

$$h_{jln}(\gamma_1)^{-1}h_{jl}(\gamma_1)h_{ik}(\gamma_0)^{-1}g_{ik}(\gamma_0)h_{ln}(\gamma_1)h_{km}(\gamma_0)^{-1}g_{ik}(\gamma_0)^{-1}h_{ikm}(\gamma_0). \tag{5.43}$$

Now for any $h \in H$ and any $g_1 \in G$ for which $\tau(g_1) = h_1$ we have

$$g_1hg_1^{-1} = \alpha\big(\tau(h_1)\big)(h) = h_1hh_1^{-1} \tag{5.44}$$

by the crossed module relation (2.5).

Using this we have

$$\begin{aligned}
&h_{jln}(\gamma_1)^{-1}h_{jl}(\gamma_1)h_{ik}(\gamma_0)^{-1}g_{ik}(\gamma_0)h_{ln}(\gamma_1)h_{km}(\gamma_0)^{-1}g_{ik}(\gamma_0)^{-1}h_{ikm}(\gamma_0)\\
&= h_{jln}(\gamma_1)^{-1}h_{jl}(\gamma_1)\ \ h_{ik}(\gamma_0)^{-1}h_{ik}(\gamma_0)\ \ h_{ln}(\gamma_1)\ h_{km}(\gamma_0)^{-1}h_{ik}\cdot\\
&\qquad\qquad\qquad\qquad \cdot(\gamma_0)^{-1}h_{ikm}(\gamma_0)\\
&= h_{jln}(\gamma_1)^{-1}h_{jl}(\gamma_1)h_{ln}(\gamma_1)\quad h_{km}(\gamma_0)^{-1}h_{ik}(\gamma_0)^{-1}h_{ikm}(\gamma_0)\\
&= h_{jn}(\gamma_1)\ \ \big(h_{ikm}(\gamma_1)^{-1}h_{ik}(\gamma_0)h_{km}(\gamma_0)\big)^{-1} \qquad \text{using (5.26 )}\\
&= h_{jn}(\gamma_1)h_{im}(\gamma_0)^{-1},
\end{aligned} \tag{5.45}$$

which we recognize to be the $H$-component of

$$\theta^{j\,|l|\,n}_{i\,|k|\,m}(\gamma) = \Big(h_{jn}(\gamma_1)h_{im}(\gamma_0)^{-1}, g_{im}(\gamma_0)\Big). \tag{5.46}$$

Thus the function $h_{\mathbb{T}}$ is indeed the function which implements the natural transformation as described in Proposition 3.3. □

5.3. **Functorial cocycles in the large.** The *functorial cocycle* given by the system of functors $\{\theta^{jn}_{ik}\}$ is a categorical analog of the traditional $G$-valued cocycle $\{g_{ij}\}_{i,j\in I}$, a family of $G$-valued (suitably regular) functions $U_i \longrightarrow G$, with $i$ running over an indexing set and $\{U_i\}_{i\in I}$ an open covering of a space $B$. In this work we will not pursue the task of constructing the global categorical structure as a counterpart to the traditional principal bundle that arises from a traditional $G$-valued cocycle.

## 6. Twisted-Product Categorical Bundles

Consider a categorical principal $\mathbf{G}$-bundle $\mathbf{P} \longrightarrow \mathbf{U}$ that is a topological principal bundle both at the object and at the morphism levels. Suppose, moreover, that both the object bundle and the morphism bundle are trivial. Nonetheless, the composition of morphisms may still be nontrivial in the sense that $\mathbf{P} \longrightarrow \mathbf{U}$ is not isomorphic to the categorical bundle $\mathbf{U}\times\mathbf{G} \longrightarrow \mathbf{U}$. A product category has a very special structure and a general categorical bundle need not be locally isomorphic, in a functorial way, to a product categorical bundle. For this reason we explore now a different type of local structure for categorical bundles.

6.1. **Twisted triviality.** Consider a categorical principal bundle

$$\pi : \mathbf{P} \longrightarrow \mathbf{U}$$

with structure categorical Lie group $\mathbf{G}$. Suppose that the object bundle $P \longrightarrow U$, where $P = \mathrm{Obj}(\mathbf{P})$ and $U = \mathrm{Obj}(\mathbf{U})$, is a smooth principal bundle. As usual, there is the associated Lie crossed module $(G, H, \alpha, \tau)$. The base category $\mathbf{U}$ has as morphisms paths on $U$. Suppose now that $P \longrightarrow U$ is isomorphic to the product bundle $U \times G \longrightarrow U$ by means of a $G$-equivariant diffeomorphism

$$F : U \times G \longrightarrow P. \tag{6.1}$$

Furthermore, we assume that there is a bijection

$$F : \mathrm{Mor}(\mathbf{U}) \times \mathrm{Mor}(\mathbf{G}) \longrightarrow \mathrm{Mor}(\mathbf{P}) \tag{6.2}$$

(we use the same notation $F$ as before) which is equivariant in the sense that

$$F(\gamma, hg)h_1 g_1 = F(\gamma, hgh_1 g_1) \tag{6.3}$$

and satisfies

$$\pi \circ F(\gamma, hg) = \gamma$$

for all $\gamma \in \mathrm{Mor}(\mathbf{U})$ and $h, h_1 \in H$ and $g, g_1 \in G$. The latter condition means that $F$ and the projection maps $\pi$ satisfy

$$\pi \circ F = \pi, \tag{6.4}$$

where $\pi$ on the left is from $\mathbf{P} \longrightarrow \mathbf{U}$ and on the right it is $\mathbf{U} \times \mathbf{G} \longrightarrow \mathbf{U}$. There is still considerable freedom in the choice of $F$; for example, we could replace it by $F_1$ given by

$$F_1(\gamma, hg) = F(\gamma, x_\gamma)hg$$

where $x_\gamma$ is some arbitrary element in $H \rtimes_\alpha G$. To fix this degree of freedom we impose the further restriction that $F$ maps the usual source of $(\gamma, e)$ to the source of $F(\gamma, e)$:

$$sF(\gamma, e) = F\big(s(\gamma), e\big) \tag{6.5}$$

for all $\gamma \in \mathrm{Mor}(\mathbf{U})$. In view of the equivariance of $F$ this implies

$$\begin{aligned} sF(\gamma, hg) = s\big(F(\gamma, e)hg\big) &= s\big(F(\gamma, e)\big)s(hg) \\ &= F\big(s(\gamma, e)\big)g \\ &= F\big(s(\gamma, e)g\big), \end{aligned} \tag{6.6}$$

wherein the second equality follows from the functorial nature of the right action of $\mathbf{G}$ on $\mathbf{P}$:

$$\mathbf{P} \times \mathbf{G} \longrightarrow \mathbf{P}.$$

When all these conditions are satisfied we will say that $F$ provides a *trivialization* for the categorical principal bundle $\pi : \mathbf{P} \longrightarrow \mathbf{U}$ and we say this bundle is *twisted-trivial*.

Let us note a crucial distinction between trivial categorical bundles and those that are twisted-trivial. There need not be a functorial isomorphism between a twisted-trivial categorical bundle and a categorical product bundle. Thus different twisted-trivial categorical bundles over the same category with the same categorical structure group need not be categorically isomorphic. The point is that the mapping $F$ discussed above is not required to be such that the target of $F(\gamma, hg)$ is $F\big(t(\gamma), t(hg)\big)$.

6.2. **Twisted-product categorical principal bundle.** Using the functions $F$ from (6.1) and (6.2) we form a category

$$\mathbf{U} \times_F \mathbf{G},$$

with object set and morphism set given by

$$\begin{aligned} \mathrm{Obj}(\mathbf{U} \times_F \mathbf{G}) &= \mathrm{Obj}(\mathbf{U}) \times \mathrm{Obj}(\mathbf{G}) = U \times G \\ \mathrm{Mor}(\mathbf{U} \times_F \mathbf{G}) &= \mathrm{Mor}(\mathbf{U}) \times \mathrm{Mor}(\mathbf{G}) = \mathrm{Mor}(\mathbf{U}) \times (H \times G), \end{aligned} \tag{6.7}$$

and source/target maps and composition defined by requiring that $F$ be a functor

$$\mathbf{U} \times_F \mathbf{G} \longrightarrow \mathbf{P}.$$

(Let us note again that $F$ is not meant to be a functor from the product category $\mathbf{U} \times \mathbf{G}$ to $\mathbf{P}$.) In more detail, this means that the source and target maps for $\mathbf{U} \times_F \mathbf{G}$ are defined by

$$\begin{aligned} s_F(\gamma, hg) &\stackrel{\text{def}}{=} F^{-1} sF(\gamma, hg) = \big(s(\gamma), g\big) \qquad \text{(using (6.6))} \\ t_F(\gamma, hg) &\stackrel{\text{def}}{=} F^{-1} tF(\gamma, hg) = F^{-1} t\big(F(\gamma, e)hg\big) \\ &= F^{-1} tF(\gamma, e)\tau(h)g. \end{aligned} \tag{6.8}$$

Since $\pi : \mathbf{P} \longrightarrow \mathbf{U}$ is a functor it commutes with the target maps:

$$\pi t F(\gamma, e) = t\pi F(\gamma, e) = t(\gamma), \tag{6.9}$$

where we also used the condition $\pi \circ F = \pi$ noted in (6.4). From this we have

$$\pi F^{-1} tF(\gamma, e) = \pi tF(\gamma, e) = t(\gamma). \tag{6.10}$$

Thus the object $F^{-1}tF(\gamma, e)$ in $\mathbf{U} \times_F \mathbf{G}$ is of the form $\big(t(\gamma), \eta_F(\gamma)\big)$ for some element

$$\eta_F(\gamma) \in G.$$

In other words,

$$F^{-1} tF(\gamma, e) = \big(t(\gamma), \eta_F(\gamma)\big) \in U \times G. \tag{6.11}$$

Consequently

$$t_F(\gamma, hg) = \big(t(\gamma), \eta_F(\gamma)\tau(h)g\big). \tag{6.12}$$

Let us now understand composition of morphisms in $\mathbf{U}\times_F\mathbf{G}$. Consider morphisms

$$(\gamma_1, e) \qquad \text{and} \qquad (\gamma_2, g),$$

where we have source-target matching:

$$\eta_F(\gamma_1) = g.$$

The projection of the composite

$$(\gamma_2, g)\circ_{\eta_F}(\gamma_1, e)$$

onto $\mathrm{Mor}(\mathbf{U})$ should be $\gamma_2 \circ \gamma_1$, and the source of the composite should be $s_F(\gamma_1, e) = \big(s(\gamma_1), e\big)$. Thus the composite is of the form

$$(\gamma_2 \circ \gamma_1, h'),$$

where $h'$ is some element of $H$ determined by $\gamma_1$ and $\gamma_2$. The target of the composite should be the same as the target of $(\gamma_2, g)$ and so

$$\eta_F(\gamma_2 \circ \gamma_1)\tau(h') = \eta_F(\gamma_2)g = \eta_F(\gamma_2)\eta_F(\gamma_1). \tag{6.13}$$

Thus this is a restriction on the element $h'$.

More general composites can then obtained by using the requirement that the usual right action

$$(\mathbf{U} \times_F \mathbf{G}) \times \mathbf{G} \longrightarrow \mathbf{U} \times_F \mathbf{G}, \tag{6.14}$$

is functorial. Thus the composite

$$(\gamma_2, h_2g_2)\circ_{\eta_F}(\gamma_1, h_1g_1) \tag{6.15}$$

is meaningful if the source-target match

$$\eta_F(\gamma_1)\tau(h_1)g_1 = g_2 \tag{6.16}$$

holds, and then the composite can be computed as follows, writing $g = \eta_F(\gamma_1)$,

$$\begin{aligned}(\gamma_2, h_2g_2)\circ_{\eta_F}(\gamma_1, h_1g_1) &= (\gamma_2, g)g^{-1}h_2g_2\circ_{\eta_F}(\gamma_1, e)h_1g_1\\ &= \Big((\gamma_2, g)\circ_{\eta_F}(\gamma_1, e)\Big)\big((g^{-1}h_2g\,g^{-1}g_2) \circ h_1g_1\big)\\ &= \Big((\gamma_2, g)\circ_{\eta_F}(\gamma_1, e)\Big)g^{-1}h_2gh_1g_1\\ &\text{(using the composition law (3.4) in } \mathrm{Mor}(\mathbf{G}))\\ &= (\gamma_2 \circ \gamma_1, h'g^{-1}h_2gh_1g_1)\end{aligned} \tag{6.17}$$

The element $h' \in H$ is determined by $\gamma_1$ and $\gamma_2$. The associative law for composition of morphisms puts a restriction on the behavior of this element.

Let us focus on the case where $h' = e$ for all $\gamma_1$ and $\gamma_2$. In view of (6.13) this means that the mapping

$$\eta_F : \mathrm{Mor}(\mathbf{U}) \longrightarrow G$$

is a 'homomorphism' in the sense that it carries composition of morphisms to the composition law in $G$.

**Proposition 6.1.** *Let* $\mathbf{U}$ *be a category,* $\mathbf{G}$ *a categorical group with associated crossed module* $(G, H, \alpha, \tau)$*. Let*

$$\eta : \mathrm{Mor}(\mathbf{U}) \longrightarrow G = \mathrm{Obj}(\mathbf{G})$$

*be a map satisfying*

$$\eta(\gamma_2 \circ \gamma_1) = \eta(\gamma_2)\eta(\gamma_1) \tag{6.18}$$

*for all* $\gamma_1, \gamma_2 \in \mathrm{Mor}(\mathbf{U})$*. Then there is a category* $\mathbf{U} \times_\eta \mathbf{G}$ *whose object set is* $\mathrm{Obj}(\mathbf{U}) \times \mathrm{Obj}(\mathbf{G})$*, with morphism set* $\mathrm{Mor}(\mathbf{U}) \times \mathrm{Mor}(\mathbf{G})$*, source and target maps given by*

$$\begin{aligned} s_\eta(\gamma, hg) &= \big(s(\gamma), g\big) \\ t_\eta(\gamma, hg) &= \big(t(\gamma), \eta(\gamma)\tau(h)g\big), \end{aligned} \tag{6.19}$$

*and composition by*

$$(\gamma_2, h_2g_2) \circ_\eta (\gamma_1, h_1g_1) = (\gamma_2 \circ \gamma_1, g^{-1}h_2gh_1g_1), \tag{6.20}$$

*where* $g = \eta_F(\gamma_1)$*, and the source* $s_\eta(\gamma_2, h_2g_2)$ *is the target* $t_\eta(\gamma_1, h_1g_1)$*. Moreover, the projection on the first factor*

$$\pi_\eta : \mathbf{U} \times_\eta \mathbf{G} \longrightarrow \mathbf{U},$$

*along with the usual right action of* $\mathbf{G}$ *on* $\mathbf{U} \times_\eta \mathbf{G}$*, makes* $\pi_F$ *a categorical principal bundle.*

Let us explore the $\eta$-twisted product $\mathbf{U} \times_\eta \mathbf{G}$ of Proposition 6.1 a bit further before looking at the example of main interest. There is an 'action' of the category $\mathbf{U}$ on the category $\mathbf{G}$ given by means of the mapping

$$E_\eta : \mathrm{Mor}(\mathbf{G}) \times \mathrm{Mor}(\mathbf{U}) \longrightarrow \mathrm{Mor}(\mathbf{G}) : (\phi, \gamma) \mapsto 1_{\eta(\gamma)^{-1}}\phi, \tag{6.21}$$

where on the right we have the multiplication in the group $\mathrm{Mor}(\mathbf{G})$ and $1_g$ denoted the identity morphism at any object $g \in \mathrm{Obj}(\mathbf{G})$. This is an action in the sense that it has the following properties:

(i) Behavior of the identity in $\mathbf{U}$:

$$E_\eta\big(\phi, \mathrm{id}_a\big) = \phi,$$

which holds because the homomorphism property of $\eta$ implies that $\eta(\mathrm{id}_a)$ is the identity in the group $\mathrm{Mor}(\mathbf{G})$ for every object $a$ in $\mathbf{U}$.

(ii) Behavior of the identity in $\mathbf{G}$:

$$E_\eta(1_g, \gamma) = 1_{\eta(\gamma)^{-1}g}$$

(iii) Behavior of composition in $\mathbf{U}$:

$$E_\eta(\phi, \gamma_2 \circ \gamma_1) = E_\eta\big(E_\eta(\phi, \gamma_2), \gamma_1\big)$$

(iv) behavior of composition in $\mathbf{U}$:

$$E_\eta\big(\phi_2 \circ \phi_1, \gamma) = E_\eta(\phi_2, \gamma) \circ E_\eta(\phi_1, \gamma),$$

which we can verify using the interchange law:

$$1_g(\phi_2 \circ \phi_1) = (1_g \circ 1_g)(\phi_2 \circ \phi_1) = (1_g\phi_2) \circ (1_g\phi_1).$$

The composition law in $\mathbf{U} \times_\eta \mathbf{G}$ is given by

$$(\gamma_2, \phi_2) \circ_\eta (\gamma_1, \phi_1) = \big(\gamma_2 \circ \gamma_1, E_\eta(\phi_2, \gamma_1) \circ \phi_1\big) \tag{6.22}$$

This structure is a slightly more general form of the *twisted semi-direct product* introduced in [14, sec. 6]

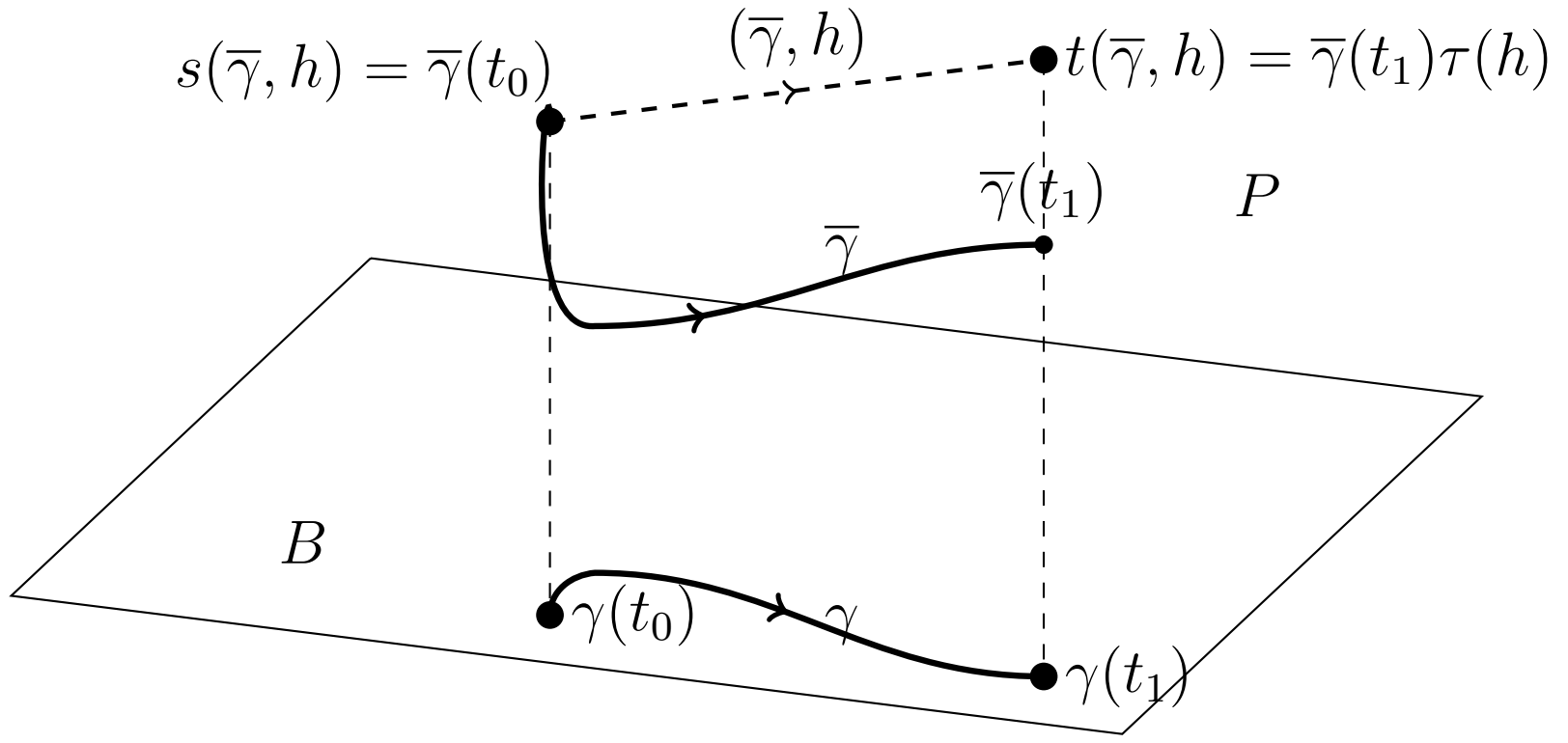

FIGURE 4. The decorated bundle, showing a morphism $(\overline{\gamma}, h)$ with its source and target

6.3. **Decorated bundles.** Let us review the notion of a *decorated bundle* from our earlier work [12]. Consider a principal $G$-bundle $\pi : P \longrightarrow B$ with a connection $\overline{A}$. Let $\mathcal{P}_{\bar{A}}P$ be the set of all $\overline{A}$-horizontal paths on $P$ and $\mathcal{P}(B)$ the set of all paths on $B$. (For expository convenience we omit technical details in this summary description; very briefly, the paths here are $C^\infty$ maps with domain $[t_0, t_1]$, constant near $t_0$ and $t_1$, and two paths that differ by a translation of the domain $t \mapsto t + k$ are identified.) Thus there is the projection map

$$\pi_{\overline{A}} : \mathcal{P}_{\bar{A}}P \longrightarrow \mathcal{P}(B) : \overline{\gamma} \mapsto \gamma = \pi \circ \overline{\gamma},$$

and the group $G$ has a natural right action on $\mathcal{P}_{\bar{A}}P$ by by right translations: $(\overline{\gamma}, g) \mapsto \overline{\gamma}g$. In this way $\mathcal{P}_{\bar{A}}P$ is a principal $G$-bundle over $\mathcal{P}M$. From this we construct a categorical principal bundle

$$\mathbf{P}^d \longrightarrow \mathbf{B}, \tag{6.23}$$

the *decorated bundle*, with structure categorical group $\mathbf{G}$; this is illustrated in Figure 4. At the object level this is just the principal $G$-bundle $\pi : P \longrightarrow B$. At the morphism level, $\mathrm{Mor}(\mathbf{B}) = \mathcal{P}(B)$ and $\mathrm{Mor}(\mathbf{P}^d) = \mathcal{P}_{\bar{A}}P \times H$, a typical element being of the form $(\overline{\gamma}, h)$ where $\overline{\gamma}$ is an $\overline{A}$-horizontal path on $P$ and $h \in H$. The action of $\mathrm{Mor}(\mathbf{G}) = H \rtimes_\alpha G$ is given by

$$(\overline{\gamma}, h)h_1g_1 = (\overline{\gamma}g_1, g_1^{-1}hh_1g_1), \tag{6.24}$$

with obvious notation. The projection of $(\overline{\gamma}, h)$ onto $\mathrm{Mor}(\mathbf{B})$ is just $\pi \circ \overline{\gamma}$. Source and target are given by

$$\begin{aligned} s(\overline{\gamma}, h) &= s(\overline{\gamma}) = \overline{\gamma}(t_0) \\ t(\overline{\gamma}, h) &= \overline{\gamma}(t_1)\tau(h) \end{aligned} \tag{6.25}$$

and composition is given by

$$(\overline{\gamma}_2, h_2) \circ (\overline{\gamma}_1, h_1) = (\overline{\gamma}_3, h_1h_2), \tag{6.26}$$

whenever $t(\overline{\gamma}_1, h_1) = s(\overline{\gamma}_2, h_2)$, where $\overline{\gamma}_3$ is the composite of $\overline{\gamma}_1$ with the right translate $\overline{\gamma}_2\tau(h_1)^{-1}$:

$$\overline{\gamma}_3 = \overline{\gamma}_2\tau(h_1)^{-1} \circ \overline{\gamma}_1. \tag{6.27}$$

For more we refer to [12] and (with notation as here) [13].

6.4. **Decoration of product bundles.** We focus now on the case when

$$P = B \times G,$$

as a product bundle over $B$. A connection form can then be specified by an $L(G)$-valued 1-form $A_0$. The horizontal lift of a path $\gamma : [t_0, t_1] \longrightarrow B$, with initial point $\big(\gamma(t_0), g\big) \in P$ is the path

$$\overline{\gamma}_g : [t_0, t_1] \longrightarrow P : u \mapsto \overline{\gamma}_g(u) \stackrel{\text{def}}{=} \big(\gamma(u), g_\gamma(u)g\big), \tag{6.28}$$

where $[t_0, t_1] \longrightarrow G : u \mapsto g_\gamma(u)$ is the solution of the differential equation

$$g'_\gamma(u)g_\gamma(u)^{-1} = -A_0\big(\gamma'(u)\big), \tag{6.29}$$

with $g_\gamma(t_0) = e$. We will denote by $\tilde{\gamma}$ the specific horizontal lift of $\gamma$ with initial point $\big(\gamma(t_0), e\big)$; thus

$$\tilde{\gamma}(u) = \big(\gamma(u), g_\gamma(u)\big). \tag{6.30}$$

6.5. **Decorated products as twisted products.** We continue to work with the preceding framework. Let us construct a mapping $\eta : \mathrm{Mor}(\mathbf{B}) \longrightarrow G$ satisfying (6.18). Let $\eta(\gamma)$ be parallel-transport along $\gamma$ by $A_0$. In more detail, for a path

$$\gamma : [t_0, t_1] \longrightarrow B,$$

we take

$$\eta(\gamma) = g_\gamma(t_1).$$

An object of $\mathbf{B} \times_F \mathbf{G}$ is of the form $(b, g) \in B \times G$, and a morphism is of the form

$$(\gamma, hg) \in \mathrm{Mor}(\mathbf{B}) \times (H \rtimes_\alpha G).$$

Now let

$$\overline{\gamma} : [t_0, t_1] \longrightarrow B \times G : t \mapsto \big(\gamma(t), g_\gamma(t)\big)$$

be the horizontal lift, initiating at $\big(\gamma(t_0), e\big)$, of $\gamma$ to the principal bundle $B \times G$ by means of the connection specified by $A_0$. Consider the bijection

$$\Theta : \mathrm{Mor}(\mathbf{P}^d) \longrightarrow \mathrm{Mor}(\mathbf{B} \times_\eta \mathbf{G}) = \mathrm{Mor}(\mathbf{B}) \times \mathrm{Mor}(\mathbf{G}) : (\overline{\gamma} g, h) \mapsto (\gamma, gh). \tag{6.31}$$

We can check that this preserves sources and targets:

$$\begin{aligned} s_\eta(\gamma, gh) &= \big(\gamma(t_0), g\big) = s(\overline{\gamma} g, h) \\ t_\eta(\gamma, gh) &= \big(\gamma(t_1), g_\gamma(t_1) t(gh)\big) = \big(\gamma(t_1), g_\gamma(t_1) g\tau(h)\big) = t(\overline{\gamma} g)\tau(h) \end{aligned} \tag{6.32}$$

It carries identity morphisms to identity morphisms. Let us now verify that it also preserves composition of morphisms. We will need to use the composition law (6.20) with $gh$ in place of $hg$; noting that $gh = ghg^{-1} \cdot g$, we have from (6.20):

$$\begin{aligned} (\gamma_2, g_2h_2) \circ_\eta (\gamma_1, g_1h_1) = (\gamma_2 \circ \gamma_1, \eta(\gamma_1)^{-1} g_2h_2g_2^{-1} \eta(\gamma_1) g_1h_1g_1^{-1}\, g_1) \\ \text{(using (6.20))} \end{aligned} \tag{6.33}$$

where, to ensure source-target matching,

$$\big(s(\gamma_2), g_2\big) = \big(t(\gamma_1), g_\gamma(t_1) g_1\tau(h_1)\big) = \big(t(\gamma_1), \eta(\gamma) g_1\tau(h_1)\big).$$

Inserting this relation in the second component of the right side in (6.33) we obtain for the second component:

$$\begin{aligned}\eta(\gamma_1)^{-1}\cdot\eta(\gamma_1)g_1\tau(h_1)h_2\big(\eta(\gamma_1)g_1\tau(h_1)\big)^{-1}\eta(\gamma_1)g_1h_1 \\ &= g_1\tau(h_1)h_2\tau(h_1)^{-1}g_1^{-1}g_1h_1\\ &= g_1h_1h_2h_1^{-1}h_1\\ &= g_1h_1h_2.\end{aligned} \tag{6.34}$$

Thus

$$(\gamma_2, g_2h_2)\circ_\eta(\gamma_1, g_1h_1) = (\gamma_2\circ\gamma_1, g_1h_1h_2). \tag{6.35}$$

On the other, we have

$$\begin{aligned}(\overline{\gamma}_2g_2, h_2)\circ(\overline{\gamma}_1g_1, h_1) &= \big((\overline{\gamma}_2g_2)\circ(\overline{\gamma}_1g_1), h_1h_2\big)\\ &\qquad\text{(using (6.26))}.\end{aligned} \tag{6.36}$$

The path $(\overline{\gamma}_2g_2)\circ(\overline{\gamma}_1g_1)$ is $\overline{A}$-horizontal, projects down to $\gamma_2\circ\gamma_1$, and has initial point $\overline{\gamma}_1(t_0)g_1 = \big(\gamma(t_0), g_1\big)$. Hence

$$(\overline{\gamma}_2g_2)\circ(\overline{\gamma}_1g_1) = \overline{\gamma_2\circ\gamma_1}g_1, \tag{6.37}$$

where $\overline{\gamma_2\circ\gamma_1}$ is the $\overline{A}$-horizontal lift of $\gamma_2\circ\gamma_1$ starting at $\big(\gamma_1(t_0), e\big)$. Then from the definition (6.31) of $\Theta$ we have

$$\begin{aligned}\Theta\Big(\big(\overline{\gamma}_2g_2, h_2\big)\circ\big(\overline{\gamma}_1g_1, h_1\big)\Big) &= \Theta\big((\overline{\gamma}_2g_2)\circ(\overline{\gamma}_1g_1),\ h_1h_2\big)\\ &= \Theta\big(\overline{\gamma_2\circ\gamma_1}g_1,\ h_1h_2\big)\qquad\text{(by (6.37))}\\ &= (\gamma_2\circ\gamma_1,\ g_1h_1h_2)\qquad\text{(using (6.31))}\\ &= (\gamma_2, g_2h_2)\circ_\eta(\gamma_1, g_1h_1)\qquad\text{(by (6.35))}\\ &= \Theta(\overline{\gamma}_2g_2, h_2)\circ_\eta\Theta(\overline{\gamma}_1g_1, h_1).\end{aligned} \tag{6.38}$$

Thus we have proved:

**Proposition 6.2.** *There is an isomorphism $\Theta$ from the decorated product categorical bundle $\mathbf{P}^d \longrightarrow \mathbf{B}$ to the twisted-product categorical bundle $\mathbf{B}\times_\eta\mathbf{G}\longrightarrow\mathbf{B}$ specified on morphisms by (6.31) and on objects by the identity map on* $\mathrm{Obj}(\mathbf{P}^d) = B\times G$.

## 7. Concluding remarks

In this paper we have (i) explored the nature of product categorical principal bundles $\mathbf{B}\times\mathbf{G}\longrightarrow\mathbf{B}$, (ii) studied functorial cocycles analogous to group-valued cocycles used as transition functions for traditional principal bundles, (iii) introduced the notion of a twisted-product categorical bundle, and (iv) proved that every categorical principal bundle arising from a traditional trivial principal bundle by decorating (by means of a connection, for example) is categorically isomorphic to a twisted-product categorical bundle. These ideas should be useful in developing a fruitful notion of local triviality for categorical principal bundles.

**Acknowledgments.** Sengupta acknowledges research support from NSA grant H98230-13-1-0210, and the SN Bose National Centre for its hospitality during visits when part of this work was done. Chatterjee acknowledges support through a fellowship from the Jacques Hadamard Mathematical Foundation. Much of Chatterjee's work towards this paper was done while he was at the Institut des Hautes Études Scientifiques, France.

## References

[1] Hossein Abbaspour, Friedrich Wagemann, *On 2-Holonomy*, (2012) `http://arxiv.org/abs/1202.2292`

[2] Paolo Aschieri, Luigi Cantini, Branislav Jurco, *Nonabelian Bundle Gerbes, their Differential Geometry and Gauge Theory*, Commun. Math. Phys., **254** (2005) 367-400.

[3] Romain Attal, *Two-dimensional parallel transport : combinatorics and functoriality*, `http://arxiv.org/abs/math-ph/0105050`

[4] Romain Attal, *Combinatorics of non-abelian gerbes with connection and curvature*, Ann. Fond. Louis de Broglie 29 (2004), no. 4, 609633.

[5] John C. Baez and Urs Schreiber, *Higher gauge theory*, Categories in algebra, geometry and mathematical physics, (eds. A. Davydov et al), Contemp. Math., Amer. Math. Soc., **431**(2007), 7-30

[6] John C. Baez and Derek K. Wise, *Teleparallel Gravity as a Higher Gauge Theory*, `http://arxiv.org/abs/1204.4339`

[7] John W. Barrett, *Holonomy and Path Structures in General Relativity and Yang-Mills Theory*, International Journal of Theoretical Physics **30** (1991) 1171-1215.

[8] Toby Bartels, *Higher Gauge Theory I: 2-Bundles*, at `http://arxiv.org/abs/math/0410328`

[9] Breen, Lawrence; Messing, William, *Differential geometry of gerbes*, Adv. Math. 198 (2005), no. 2, 732846.

[10] A. Caetano, and Roger F. Picken, *An axiomatic definition of holonomy.* Internat. J. Math. 5 (1994), no. 6, 835848.

[11] Saikat Chatterjee, Amitabha Lahiri, and Ambar N. Sengupta, *Parallel Transport over Path Spaces*, Reviews in Math. Phys. **9** (2010) 1033-1059.

[12] Saikat Chatterjee, Amitabha Lahiri, and Ambar N. Sengupta, *Path space Connections and Categorical Geometry*, Journal of Geom. Phys. **75** (2014) 1033-1059.

[13] Saikat Chatterjee, Amitabha Lahiri and Ambar N. Sengupta, *Twisted actions of categorical groups*, Theory and Applications of Categories, Vol. 29, 2014, No. 8, pp 215-255.

[14] Saikat Chatterjee, Amitabha Lahiri and Ambar N. Sengupta, *Connections on Decorated Path Bundles*, (2014).

[15] Magnus Forrester-Barker, *Group Objects and Internal Categories*, at `http://arXiv:math/0212065v1`

[16] Krzystof Gawędzki and Nuno Reis, *WZW branes and gerbes.* Rev. Math. Phys. 14 (2002), no. 12, 12811334.

[17] Gregory Maxwell Kelly and Ross Street, *Reviews of the elements of* 2-*categories*, Springer Lecture Notes in Mathematics, **420**, Berlin, (1974), 75-103.

[18] Saunders Mac Lane, *Categories for the Working Mathematician.* Springer-Verlag (1971).

[19] João Faria Martins and Roger Picken, *On two-Dimensional Holonomy*, Trans. Amer. Math. Soc. 362 (2010), no. 11, 5657-5695, `http://arxiv.org/abs/0710.4310`

[20] João Faria Martins and Roger Picken, *Surface holonomy for non-abelian 2-bundles via double groupoids*, Advances in Mathematics **226** (2011), no. 4, 3309-3366.

[21] Michael K. Murray, An introduction to bundle gerbes. *The many facets of geometry*, 237260, Oxford Univ. Press, Oxford, 2010.

[22] Arthur Parzygnat, *Gauge Invariant Surface Holonomy and Monopoles*, `http://arxiv.org/pdf/1410.6938v1.pdf`

[23] Urs Schreiber and Konrad Waldorf, *Parallel Transport and Functors*, J. Homotopy Relat. Struct. **4**, (2009), no. 1, 187-244

[24] Urs Schreiber and Konrad Waldorf, *Connections on non-abelian Gerbes and their Holonomy*, Theory and Applications of Categories **28** (2013), 476-540

[25] Emanuele Soncini and Roberto Zucchini, *A New Formulation of Higher Parallel Transport in Higher Gauge Theory*, `http://arxiv.org/abs/1410.0775`

[26] David Viennot, *Non-abelian higher gauge theory and categorical bundle*, (2012) `http://arxiv.org/abs/1202.2280`

School of Mathematics, Indian Institute of Science Education and Research, CET Campus, Thiruvananthapuram, Kerala-695016, India
*E-mail address*: `saikat.chat01@gmail.com`

Amitabha Lahiri, S. N. Bose National Centre for Basic Sciences, Block JD, Sector III, Salt Lake, Kolkata 700098, West Bengal, India
*E-mail address*: `amitabhalahiri@gmail.com`

Ambar N. Sengupta, Department of Mathematics, Louisiana State University, Baton Rouge, Louisiana 70803, USA
*E-mail address*: `ambarnsg@gmail.com`